\documentclass[12pt,reqno]{amsart}
\usepackage{amsmath}
\usepackage{amsfonts}
\usepackage{amssymb}
\usepackage{caption}
\usepackage{latexsym}
\usepackage{color}
\usepackage{hyperref}
\usepackage{url}
\usepackage{enumerate}
\usepackage{booktabs}
\usepackage{longtable}
\usepackage{multirow}
\usepackage{lscape}
\usepackage{tabularx}
\usepackage[figuresleft]{rotating}
\usepackage{adjustbox}

\usepackage{xcolor} 
\usepackage{colortbl}

\renewcommand{\proof}{\noindent{\it Proof.\ \ }}
\renewcommand{\qed}{\ifmmode\square\else\nolinebreak\hfill
$\Box$\fi\par\vskip12pt}

 \renewcommand\a{\alpha}  \renewcommand\b{\beta}  
 \newcommand\A{\mathrm{A}} \newcommand\AGL{\mathrm{AGL}}   \newcommand\ASL{\mathrm{ASL}}  \newcommand\Aut{\mathrm{Aut}}
 
\newcommand\C{\mathbf{C}}   \newcommand\Cay{\mathrm{Cay}}    
\newcommand\D{\mathrm{D}}

\newcommand\F{\mathrm{F}}  
\def\G{\mathrm{G}} \newcommand\GammaL{\mathrm{\Gamma L}}    \newcommand\GL{\mathrm{GL}}  
\newcommand\Ga{\mathrm{\Gamma}}   
 
\newcommand\J{\mathrm{J}}
\newcommand\K{\mathsf{K}}
\newcommand\M{\mathrm{M}}  

 \renewcommand\O{\mathrm{O}} 
\renewcommand\P{\mathrm{P}}  \newcommand\PGL{\mathrm{PGL}} \newcommand\PGammaL{\mathrm{P\Gamma L}}     \newcommand\PSL{\mathrm{PSL}} \newcommand\PSigmaL{\mathrm{P\Sigma L}} \newcommand\PSigmaU{\mathrm{P\Sigma U}}  \newcommand\PSp{\mathrm{PSp}} 
  
\newcommand\Q{\mathrm{Q}}

 \newcommand\SL{\mathrm{SL}}   \newcommand\Sp{\mathrm{Sp}}    \newcommand\Sy{\mathrm{S}}

\newcommand\ZZ{\mathbb{Z}}

  \newcommand\AG{\mathrm{AG}}
\newcommand\diam{\mathrm{diam}}
\newcommand\BB{\mathcal{B}}

\newcommand\GG{\mathcal{G}}
\newcommand\ov{\overline}
\newcommand\la{\langle} \newcommand\ra{\rangle}

\newtheorem{theorem}{Theorem}[section]%
\newtheorem{lemma}[theorem]{Lemma}%
\newtheorem{corollary}[theorem]{Corollary}%
\newtheorem{proposition}[theorem]{Proposition}%
\newtheorem{definition}[theorem]{Definition}%
\newtheorem{example}[theorem]{Example}%
\newtheorem{question}[theorem]{Question}%

\begin{document}

\title[Geodesic transitive graphs of small valency]
{Geodesic transitive graphs of small valency}

\thanks{2020 MR Subject Classification 20B15, 20B30, 05C25.}

\author[J.-J. Huang]{Jun-Jie Huang}


\address{Jun-Jie Huang\\
School of Mathematical Sciences, Laboratory of Mathematics and Complex Systems, MOE\\
Beijing Normal University\\
Beijing \\
100875, P. R. China}
{\email{jjhuang@bnu.edu.cn, 20118006@bjtu.edu.cn (J.-J. Huang)}}
\maketitle

\begin{abstract}
  For a graph $\Gamma$, the {\em distance} $d_\Gamma(u,v)$ between two distinct vertices $u$ and $v$ in $\Gamma$ is  defined as the length of the shortest path from $u$ to $v$, and the {\em diameter} $\diam(\Gamma)$ of $\Gamma$ is the maximum distance between $u$ and $v$ for all vertices $u$ and $v$ in the vertex set of $\Gamma$. For a positive integer $s$, a path $(u_0,u_1,\ldots,u_{s})$ is called an {\em $s$-geodesic} if the distance of $u_0$ and $u_s$ is $s$. The graph $\Gamma$ is said to be {\em distance transitive} if for any vertices $u,v,x,y$ of $\Ga$ such that $d_\Ga(u,v)=d_\Ga(x,y)$, there exists an automorphism of $\Gamma$ that maps the pair $(u,v)$ to the pair $(x,y)$. Moreover, $\Gamma$ is said to be {\em geodesic transitive} if for each $i\leq \diam(\Ga)$, the full automorphism group acts transitively on the set of all $i$-geodesics. In the monograph [Distance-Regular Graphs, Section 7.5], the authors listed all distance transitive graphs of valency at most $13$. By using this classification, in this paper, we provide a complete classification of geodesic transitive graphs with valency at most $13$. As a result, there are exactly seven graphs of valency at most $13$ that are distance transitive but not geodesic transitive.
\end{abstract}

\qquad {\textsc k}{\scriptsize \textsc {eywords.}  geodesic transitive, distance transitive, automorphism group} {\footnotesize}

\section{Introduction}

All graphs considered in this paper are finite, connected, undirected and simple.
For a graph $\Ga$, we denote its vertex set, edge set and the full automorphism group by $V(\Ga)$, $E(\Ga)$ and $\Aut(\Ga)$ respectively.
The {\em distance} between two distinct vertices $u$ and $v$ in $\Ga$ is the length of the shortest path from $u$ to $v$, denoted by $d_\Ga(u,v)$,
and the {\em diameter} $\diam(\Ga)$ of $\Ga$ is the maximum value of the distance between any two vertices $u$ and $v$ in $V(\Ga)$.
For a positive integer $s$, an {\em $s$-arc} of $\Ga$ is a sequence of vertices $(u_0,u_1,\ldots,u_s)$ in $\Ga$ such that $u_i$ is adjacent to $u_{i+1}$ for $0\leq i\leq s-1$ and $u_{j-1}\neq u_{j+1}$ for $1\leq j\leq s-1$.
This $s$-arc is called an {\em $s$-geodesic} if $d_\Ga(u_0,u_s)=s$.

Let $G$ be a subgroup of $\Aut(\Ga)$ that acts transitively on $V(\Ga)$.
The graph $\Ga$ is said to be {\em $(G,s)$-arc transitive} if $\Ga$ has at least one $s$-arc and $G$ is transitive on the set of $s$-arcs of $\Ga$.
Similarly, if for each $1\leq i\leq s$, $\Ga$ has at least one $i$-geodesic, and $G$ is transitive on the set of all $i$-geodesics of $\Ga$, then $\Ga$ is called {\em $(G,s)$-geodesic transitive}.
In particular, a $(\Aut(\Ga),s)$-arc transitive or $(\Aut(\Ga),s)$-geodesic transitive graph will be simply called {\em $s$-arc transitive} or {\em $s$-geodesic transitive}, respectively.
Moreover, a $\diam(\Ga)$-geodesic transitive graph $\Ga$ is called {\em geodesic transitive},
and an $s$-arc transitive but not $(s+1)$-arc transitive graph is called {\em $s$-transitive}.
It should be noted that every $s$-geodesic of $\Ga$ is an $s$-arc, but not all $s$-arcs are $s$-geodesics.
For example, all $s$-arcs that lie in a cycle of length $2s-1$ are not $s$-geodesics.
Therefore, the family of $s$-arc transitive graphs is properly included in the family of $s$-geodesic transitive graphs.
We would like to emphasize, however, that there exist graphs that are $s$-geodesic transitive but not $s$-arc transitive,
see~\cite{HFZY,JDLP} for example.

The investigation of $s$-arc transitive graphs was initiated by Tutte~\cite{Tutte1947} and Weiss~\cite{Weiss81}.
They proved the nonexistence of $8$-arc transitive graphs of valency at least three.
Since their foundational work, $s$-arc transitive graphs have captured significant interest in algebraic graph theory,
see~\cite{CP83,DMM,FZL,Li2001,Praeger93,Zhou} and the references therein.
Consequently, when delving into $s$-geodesic transitive graphs,
the main focus is on studying those graphs that are not $s$-arc transitive.
For example, the local structure of $2$-geodesic transitive graphs was investigated in~\cite{DJLP13-2},
and Huang et al.~\cite{HFZY} classified $2$-geodesic transitive graphs of order $p^n$, where $p$ is a prime and $n \leq 3$.
Moreover, a fascinating reduction theorem for $3$-geodesic transitive but not $3$-arc transitive graphs was developed using the normal quotient strategy in~\cite{JP}.
For more results, we refer to~\cite{DJLP13,DJLP14,JDLP} for example.

More recently, $s$-geodesic transitive but not $s$-arc transitive graphs under specific valency have become a focus of research.
For $s=2$, such graphs of valencies $4$, $6$, $9$, $p$, $2p$,
or $3p$ (where $p$ is a prime) have been studied in the literature~\cite{CJS,DJLP15,Jin2015-1,Jin2015-2,Jin2015-3,JLX}.
For $s=3$, such graphs of valency at most $5$ were classified~\cite{HFWZ,Jin15,Jin18}.
From these results, it is evident that most $s$-geodesic transitive graphs are geodesic transitive,
meaning that geodesic transitive graphs constitute a significant proportion among these graphs.
In view of the above, a natural question arises:

\begin{question}\label{ques}
Classify geodesic transitive graphs of small valency.
\end{question}

It is noted that geodesic transitive graphs are distance transitive.
To explain this, we first introduce some notations and definitions as follows.
For a positive integer $i\leq\diam(\Ga)$ and a vertex $u$ of $\Ga$, denote by $\Ga_i(u)=\{v\in V(\Ga)\mid d_\Ga(u,v)=i\}$.
For a subgroup $G$ of $\Aut(\Ga)$, we say that $\Ga$ is {\em $G$-distance transitive} if for each $i\leq\diam(\Ga)$ and any vertex $u\in V(\Ga)$,
the vertex stabilizer $G_u$ acts transitively on $\Ga_i(u)$.
Moreover, a $\Aut(\Ga)$-distance transitive is simply called {\em distance transitive}.
By definition, a $G$-geodesic transitive graph is also $G$-distance transitive, but the inverse is not true, refer to Lemma~\ref{Paley}.

Let $\Gamma$ be a connected geodesic transitive graph of valency at most $13$.
Our goal of this paper is to classify such a graph $\Gamma$.
This classification task is made easier by the fact that $\Gamma$ is distance transitive.
Notably, all distance transitive graphs of valency at most $13$ are exhaustively listed in~\cite{BCN}.
Most of the graphs listed in~\cite[Section 7.5]{BCN} belong to one of the following eleven families of graphs:
the complete graph, the complete bipartite graph, the complete bipartite minus a matching, the complete multipartite graph,
the Hamming graph, the folded $d$-cube, the Odd graph, the double cover of the Odd graph, the Johnson graph, the doubled Grassmann graph,
and the $\AG(2,q)$ minus a parallel class (these graphs will be introduced in Section~\ref{example}).
The geodesic transitivity of the first ten classes of graphs as above is known, see~\cite{FH,HFZY,JDLP,JLX} for example.
In Lemma~\ref{AG2q}, it is shown that the $\AG(2,q)$ minus a parallel class is also geodesic transitive.

As note above, classifying geodesic transitive graphs of valency at most $13$ reduces to determining whether each of the graphs listed in~\cite{BCN},
excluding the eleven families of graphs, is geodesic transitive or not.
Theoretically, based on the discussions in~\cite{Jin2015-1,JLX},
it is possible to determine the geodesic transitivity of such a graph  according to its properties.
However, this approach involves examining hundreds of graphs listed in~\cite[Section 7.5]{BCN},
making the proofs both extensive and repetitive.
Moreover, the proofs for showing that a given graph is geodesic transitive usually have similar structures, further complicating the process.
To overcome these challenges and ensure the completeness and accuracy of our results,
we employ the Magma computational software to determine the geodesic transitivity of each graph.
By taking advantage of the capabilities of Magma,
we can efficiently classify $\Gamma$ and provide a conclusive answer regarding its geodesic transitivity.

\begin{theorem}\label{geo-small}
Let $\Ga$ be a connected graph of valency $k$, where $2\leq k\leq 13$.
Then $\Ga$ is geodesic transitive if and only if $\Ga$ is isomorphic to one of the graphs listed in Tables~$\ref{val3}$--$\ref{val13}$,
with the exception of the seven graphs enumerated in Lemma~$\ref{ngeo-trans}$.
\end{theorem}

We remark that a connected geodesic transitive graph of valency at most $2$ is isomorphic to the cycle $\C_n$ of length $n\geq1$.
Combining this fact with Theorem~\ref{geo-small}, we obtain the classification of all geodesic transitive graphs of valency at most $13$.
Moreover, the four graphs listed in Examples~\ref{graph22-6}--\ref{graph68-12},
together with the Paley graph $\P(q)$ with $q\geq13$ a prime power and $q\equiv1\pmod{4}$,
and the Peisert graph Pei$(p^r)$ with $r$ even and $p$ a prime such that $p\equiv3\pmod{4}$ and $p^r\neq 9$ (see~\cite{JDLP}),
are all distance transitive but not geodesic transitive graphs discovered so far.
Consequently, a natural and interesting question is: to construct distance transitive but not geodesic transitive graphs.

The layout of this paper is as follows.
The subsequent section is divided into two main parts.
The first subsection provides the definitions and related properties of groups and graphs,
including a proof of the necessary and sufficient conditions for a graph to be $s$-geodesic transitive.
The second subsection focuses on examples of geodesic transitive graphs,
specifically proving that $\AG(2,q)$ minus a parallel class is a geodesic transitive graph.
Finally, a summary of geodesic transitive graphs of valency at most $13$ is provided in Section \ref{sum},
which includes the Magma code and Tables~$\ref{val3}$--$\ref{val13}$.

\section{Preliminaries}\label{Prel}
In this section, we give some definitions concerning groups and graphs which will be used in our analysis.
Moreover, some examples of distance graphs are constructed in the section.

For a positive integer $n$, denote by $\ZZ_n$ the cyclic group of order $n$, and by $\D_{2n}$ the dihedral group of order $2n$.
For a prime $p$ and a positive integer $r$, denote by $\ZZ_p^r$ (or simply $p^r$) the elementary abelian group of order $p^r$,
and by $[p^r]$ a finite group of order $p^r$, respectively.
Let $p^{1+2}_+$ and $p^{1+2}_-$ be the nonabelian group of order $p^3$ of exponent $p$ or $p^2$ respectively.
More precisely, we have
\begin{align*}
p^{1+2}_+&=\la a,b,c\mid a^p=b^p=c^p=1,[a,b]=c,[a,c]=[b,c]=1\ra,\\
p^{1+2}_-&=\la a,b\mid a^{p^2}=b^p=1,a^b=a^{1+p}\ra.
\end{align*}
For two groups $A$ and $B$, denote by $A\times B$ the direct product of $A$ and $B$, by $A:B$ a semidirect product of $A$ by $B$, and by $A.B$ an extension of $A$ by $B$.
Moreover, we use $p^{r+t}$ to denote a graph with structure $\ZZ_p^r.\ZZ_p^t$, and more notations of group structures, we refer to Atlas~\cite{Atlas}.
For a group $G$ and a subgroup $H$ of $G$, we use $[G:H]$ denote the set of all right cosets of $H$ in $G$.

\subsection{Graph theoretical results}
We now introduce some parameters that play an important role in the study of distance transitive graphs, especially those with small valencies.

\begin{definition}\label{def}
Assume that $\Ga$ is an $s$-geodesic transitive graph with $s\leq\diam(\Ga)$.
Let $u\in V(\Ga)$, and let $v\in\Ga_i(u)$ with $i\leq \diam(\Ga)$.
Then, the number of edges from $v$ to $\Ga_{i-1}(u)$, $\Ga_i(u)$, and $\Ga_{i+1}(u)$ does not depend on the choice of $v$,
and these numbers are denoted by $c_i$, $a_i$ and $b_i$, respectively.
\end{definition}

Clearly, for an $s$-geodesic transitive graph with diameter $d$,
the parameters $a_i$, $b_i$ and $c_i$ are well-defined, and $a_i+b_i+c_i$ is equal to the valency of $\Ga$ for all $i\leq s\leq d$.
When $s=d$, the array $\{b_0,b_1,\ldots,b_{d-1};c_1,c_2,\ldots,c_{d}\}$,
where $b_0$ is the valency of $\Ga$, is called the {\em intersection array} of $\Ga$.
Some properties of these parameters are given in~\cite{Biggs,BCN}.
The following result can be directly deduced from~\cite[Lemma 2.3]{JP}.

\begin{proposition}\label{geo-tran}
Let $s$ be a positive integer, and let $\Ga$ be a $(G,s)$-geodesic transitive graph such that $b_s\leq 1$, where $G\leq\Aut(\Ga)$. Then $\Ga$ is $G$-geodesic transitive.
\end{proposition}

\medskip
Let $G$ be a finite group and let $S\subseteq G\setminus\{1\}$ be such that $S=S^{-1}:=\{s^{-1}\mid s\in S\}$ and $\la S\ra=G$.
Then the {\em Cayley graph} $\Ga=\Cay(G,S)$ is the graph with vertex set $G$ and two vertices $x,y\in G$ are adjacent if and only if $yx^{-1}\in S$.

\medskip
Let $G$ be a transitive permutation group on a set $\Omega$.
Then $G$ induces a natural action on $\Omega\times \Omega$.
The corresponding orbits are called {\em orbitals}.
Clearly, for an orbital $O$, $O^*:=\{(\b,\a)\mid (\a,\b)\in O\}$ is also an orbital.
An orbital $O$ is said to be {\em self-paired} if $O=O^*$.
It is well-known that $\Delta:=\{(\a,\a)\mid\a\in\Omega\}$ is a self-paired orbital, called {\em diagonal}.
If $O$ is non-diagonal self-paired, then the {\em orbital graph} of $G$ corresponding to $O$ has vertex set $\Omega$ and edge set $\{\{u,v\}\ |\ (u,v)\in O\}$.
It is easy to show that this orbital graph is arc-transitive.

For each orbital $O$ of $G$ and each $\a\in O$, define $O(\a)=\{\b\in\Omega\mid(\a,\b)\in O\}$.
Then $\O(\a)$ is an orbit of $G_\a$ acts on $\Omega$, called {\em suborbital} of $G$.
When $O$ is self-paired or diagonal, $O(\a)$ is called {paired suborbit} or {\em trivial}, respectively.
Moreover, the mapping $O\mapsto O(\a)$ is a bijection from the set of orbitals of $G$ onto the set of orbits of $G_\a$.
For an arc transitive graph $\Ga$, let $A=\Aut(\Ga)$ and let $u\in V(\Ga)$.
Then $\Ga$ can be constructed as an orbital graph of the action of $G$ on the right coset $[A:A_u]$.

\medskip
The {\em complement graph} $\ov{\Ga}$ of a graph $\Ga$ is define the graph with vertex set $V(\Ga)$ and  two vertices are adjacent in $\ov{\Ga}$ if and only if they are not adjacent in  $\Ga$.
Clearly, $\Aut(\ov{\Ga})=\Aut(\Ga)$.

\medskip
A classical approach to the study of $(G,s)$-geodesic transitive graphs is to take normal quotient graphs (see~\cite{DGLP12,JP}).
Let $N$ be a nontrivial normal subgroup of $G$ such that $N$ is intransitive on $V(\Ga)$.
Let $\BB$ be the set of all orbits of $N$ acting on $V(\Ga)$.
The {\em normal quotient graph $\Ga_N$} of $\Ga$ is the graph with vertex set $\BB$, and two vertices $B,C\in\BB$ are adjacent in $\Ga_N$ if and only if there is a vertex $u\in B$ and $v\in C$ such that $\{u,v\}$ is an edge of $\Ga$.
If further, $|\Ga(u)\cap B|=0$ or $1$ for all $u\in V(\Ga)$ and $B\in \BB$, then $\Ga$ is an {\em normal cover} (or simply a {\em cover}) of $\Ga_N$.
The following result give a reduction theorem for studying $(G,s)$-geodesic transitive graphs, we refer to~\cite[Lemma 3.2]{JP} and \cite[Theorem 4.1]{Praeger93}.

\begin{proposition}\label{cover}
Let $\Ga$ be a connected $(X,s)$-geodesic transitive graph with $X\leq\Aut(\Ga)$ and $s\geq2$.
Let $N$ be a nontrivial normal subgroup of $X$ such that $N$ acting on $V(\Ga)$ has at least $3$ orbits.
Then either $\Ga_N\cong\K_m$ and $\Ga\cong\K_{m[b]}$ with $m\geq3$ and $b\geq2$, or the following holds:
\begin{enumerate}[\em (i)]
  \item $N$ is semiregular on $V(\Ga)$, $X/N\leq\Aut(\Ga_N)$, and $\Ga$ is a cover of $\Ga_N$;
  \item $\Ga_N$ is $(X/N,s')$-geodesic transitive with $s'=\min\{\diam(\Ga_N),s\}$;
  \item $\Ga$ is $(X,s)$-arc transitive if and only if $\Ga_N$ is $(X/N,s)$-arc transitive, where $2\leq s\leq 7$.
\end{enumerate}
\end{proposition}

Let $\Ga$ be a graph.
The {\em bipartite double cover} $2.\Ga$ of $\Ga$ is the bipartite graph with biparts $\{(u,0)\mid u\in V(\Ga)\}$ and $\{(u,1)\mid u\in V(\Ga)\}$
and two vertices $(u,0)$ and $(v,1)$ are adjacent in $2.\Ga$ if $u$ is adjacent to $v$ in $\Ga$.
Then $2.\Ga$ is connected if and only if $\Ga$ is connected and non-bipartite (see~\cite[Theorem 1.11.1]{BCN}).
For each $\a\in\Aut(\Ga)$ and $u\in V(\Ga)$, define $\sigma_\a: (u,i)\mapsto (u^\a,i)$ for $i=0$ or $1$.
Then $\sigma_\a\in\Aut(2.\Ga)$.
Moreover, the mapping $\tau: (u,i)\mapsto (u,1-i)$ induces an automorphism of $2.\Ga$, still denote by $\tau$.
Clearly, $\sigma_\a\tau=\tau\sigma_\alpha$.
Let $G$ be a subgroup of $\Aut(\Ga)$, and let
\begin{equation*}
X(G)=\la\sigma_\a,\tau\mid\a\in G\ra.
\end{equation*}
Then $X(G)\cong G\times\la\tau\ra$ and $X(G)\leq \Aut(2.\Ga)$.
Consider the quotient graph $(2.\Ga)_{\la\tau\ra}$, with vertex set $\{\{(u,0),(u,1)\}\mid u\in V(\Ga)\}$.
Then $(2.\Ga)_{\la\tau\ra}\cong\Ga$, and hence $2.\Ga$ is a $\la\tau\ra$-cover of $\Ga$.

\medskip
For a graph $\Ga$ with an automorphism group $G$ and an integer $s\geq1$,
it is well-known that $\Ga$ is $(G,s+1)$-arc transitive if and only if $\Ga$ is $(G,s)$-arc transitive, and for an $s$-arc $(u_0,u_1,\ldots,u_s)$,
the stabilizer $G_{u_0,u_1,\ldots,u_s}$ is transitive on the set $\Ga(v_s)\setminus\{u_{s-1}\}$.
Below, we present a similar sufficient and necessary condition for determining whether a graph is $s$-geodesic transitive.

\begin{lemma}\label{condition}
Let $\Ga$ be a graph with an automorphism group $G$ and let $s\leq\diam(\Ga)$ be an positive integer.
Then $\Ga$ is $(G,s+1)$-geodesic transitive if and only if $s+1\leq\diam(\Ga)$, $\Ga$ is $(G,s)$-arc transitive,
and for an $s$-geodesic $(u_0,u_1,\ldots,u_{s})$, the stabilizer $G_{u_0,u_1,\ldots,u_s}$ is transitive on $\Ga(u_s)\cap\Ga_{s+1}(u_0)$.
\end{lemma}

\proof Assume that $\Ga$ is $(G,s+1)$-geodesic transitive.
Then $s+1\leq\diam(\Ga)$ and $\Ga$ is $(G,s)$-geodesic transitive.
Let $(u_0,u_1,\ldots,u_s)$ be an $s$-geodesic of $\Ga$.
Then for each two vertices $u_{s+1},v_{s+1}\in\Ga(u_s)\cap\Ga_{s+1}(u_0)$,
we know that $(u_0,\ldots,u_s,u_{s+1})$ and $(u_0,\ldots,u_s,v_{s+1})$ are two $(s+1)$-geodesics of $\Ga$.
By the $(G,s+1)$-geodesic transitivity of $\Ga$, there is an element $g\in G$ such that
\begin{equation*}
(u_0,\ldots,u_s,u_{s+1})^g=(u_0,\ldots,u_s,v_{s+1}).
\end{equation*}
It follows that $g\in G_{u_0,u_1,\ldots,u_s}$ and $u_{s+1}^g=v_{s+1}$,
and hence $G_{u_0,u_1,\ldots,u_s}$ is transitive on $\Ga(u_s)\cap\Ga_{s+1}(u_0)$.

Conversely, assume that $(u_0,u_1,\ldots,u_s,u_{s+1})$ and $(v_0,v_1,\ldots,v_s,v_{s+1})$ are two $(s+1)$-geodesics of $\Ga$.
Since $\Ga$ is $(G,s)$-geodesic transitive, we have that
\begin{equation*}
(v_0,v_1,\ldots,v_s,v_{s+1})^g=(u_0,u_1,\ldots,u_s,v_{s+1}^g)
\end{equation*}
for some $g\in G$.
Now, $u_{s+1},v_{s+1}^g\in \Ga(u_s)\cap\Ga_{s+1}(u_0)$,
and since $G_{u_0,u_1,\ldots,u_s}$ is transitive on $\Ga(u_s)\cap\Ga_{s+1}(u_0)$,
we have $(v_{s+1}^g)^h=u_{s+1}$ for some $h\in G_{u_0,u_1,\ldots,u_s}$.
Therefore,
\begin{equation*}
(v_0,v_1,\ldots,v_s,v_{s+1})^{gh}=(u_0,u_1,\ldots,u_s,u_{s+1}),
\end{equation*}
implying that $\Ga$ is $(G,s+1)$-geodesic transitive, completing the proof. \qed

Based on Lemma \ref{condition} and the definition of a $(G,s)$-geodesic transitive graph,
we can immediately obtain the following corollary.

\begin{corollary}\label{condition1}
Let $\Ga$ be a graph with an automorphism group $G$ and let $s\leq\diam(\Ga)$ be an positive integer.
Then $\Ga$ is $(G,s)$-geodesic transitive if and only if $\Ga$ is $G$-arc transitive, and for each $1\leq i\leq s$,
$G_{u_0,u_1,\ldots,u_i}$, the stabilizer of the $i$-geodesic $(u_0,u_1,\ldots,u_i)$ in $G$, is transitive on $\Ga(u_i)\cap\Ga_{i+1}(u_0)$.
\end{corollary}

For an $s$-geodesic transitive graph of valency $b_0$, and for $1\leq i\leq s$, let $b_i$ be defined as in Definition~$\ref{def}$.
Then $b_i=|\Ga(u_i)\cap\Ga_{i+1}(u_0)|$.
By Corollary~\ref{condition1} and the Orbit-stabilizer Theorem (see~\cite[Theorem 1.4A~(iii)]{Dixon}), we have the subsequent result.

\begin{corollary}\label{bound}
Let $\Ga$ be a $(G,s)$-geodesic transitive graph of valency $b_0$.
Then $b_0b_1\cdots b_{s-1}$ is a divisor of $|G_{u}|$, where $u\in V(\Ga)$ and $b_1,\ldots,b_{s-1}$ are defined as in Definition~$\ref{def}$.
\end{corollary}

\subsection{Examples of geodesic transitive graphs}\label{example}
We first introduce some notation for some well-known graphs.
Denote by $\C_n$ the cycle graph of length $n$, $\K_n$ the complete graph of order $n$, $\K_{n,n}$ the complete bipartite graph of order $2n$,
$\K_{n,n}-n\K_2$ the subgraph of $\K_{n,n}$ minus a matching, and $\K_{m[b]}$ the complete multipartite graph consisting $m\geq 3$ parts of size $b\geq 2$, respectively.
Clearly, these graphs are geodesic transitive and their full automorphism groups are known.

Let $d,n\geq2$ be two integers.
The {\em Hamming graph} $H(d,n)$ is the graph with vertex set $\{(x_1,x_2,\ldots,x_d)\mid x_i=0,1,\ldots,n-1\text{~for~}1\leq i\leq d\}$,
and two vertices $(x_1,x_2,\ldots,x_d)$ and $(y_1,y_2,\ldots,y_d)$ are adjacent if $|\{i\mid x_i\neq y_i,1\leq i\leq d\}|=1$.
In particular, $H(d,2)$ is also known as the {\em $d$-cube}.
The Hamming graph has the following properties,
we refer to~\cite[Theorem 9.2.1]{BCN}, \cite[Example 3.3]{HFZY}, \cite[P. 424 (1.2)]{IP} and~\cite[Proposition 2.2]{JDLP}.

\begin{proposition}\label{Hamming}
Let $\Ga=H(d,n)$ be the Hamming graph, where $d,n\geq2$. Then the following holds.
\begin{enumerate}[\rm (1)]
  \item $\Ga$ has order $n^d$, valency $d(n-1)$, diameter $d$, and $\Aut(\Ga)\cong\Sy_n\wr\Sy_d$.
  \item $\Ga$ is a connected $2$-transitive graph if and only if $\Ga=H(d,2)$ and $d\geq3$.
  \item $\Ga$ is geodesic transitive, has intersection array given by $b_i=(d-i)(n-1)$ and $c_i=i$ for $0\leq i\leq d$.
  \item The complement graph $\ov{H(2,n)}$, where $n\geq3$, is a connected geodesic transitive but not $2$-arc transitive graph with intersection array $\{(n-1)^2,2(n-2);1,(n-1)(n-2)\}$.
\end{enumerate}
\end{proposition}

Let $d\geq3$ be an integer.
The {\em folded $d$-cube $\square_d$} is constructed from the Hamming graph $H(d,2)$ as follow.
The vertex set of $\square_d$ is $\{\{u,u'\}: u,u'\in V(H(d,2)) \mid  d_{H(d,2)}(u,u')=d\}$,
and two vertices $\{u,u'\}$ and $\{v,v'\}$ are adjacent in $\square_d$ if and only if there exists $x\in \{u,u'\}$ and $y\in \{v,v'\}$ such that $x$ and $y$ are adjacent in  $H(d,2)$.
Then $H(d,2)$ is a cover of $\square_d$.
Moreover, $\square_3\cong\K_4$, $\square_4\cong\K_{4,4}$ and $\square_d$ has diameter $[d/2]$.
By Propositions~\ref{cover} and~\ref{Hamming} we know that $\square_d$ is geodesic transitive.
Combining this with~\cite[Theorem 7.5.2]{BCN} and \cite[P. 424 (1.2)]{IP} we obtain the next result.

\begin{proposition}\label{folded}
Let $\Ga=\square_d$ be the folded $d$-cube with $d\geq5$. Then the following holds.
\begin{enumerate}[\rm (1)]
  \item $\Ga$ has order $2^{d-1}$, valency $d$, diameter $k:=[d/2]$,
        and intersection array $\{d,d-1,\ldots,d+1-k;1,2,\ldots,d-1,c_k\}$ with $c_k=k$ if $d$ is odd and $c_k=d$ if $d$ is even.
  \item $\Ga$ is a connected $2$-transitive and geodesic transitive graph, and $\Aut(\Ga)\cong\ZZ_2^{k-1}\rtimes\Sy_k$.
\end{enumerate}
\end{proposition}

Let $\Omega$ be a finite set with cardinality $n$, and let $\Omega^{\{k\}}$ be the set of all $k$-subsets of $\Omega$.
The {\em Johnson graph $\J(n,k)$} is the graph whose vertex set $\Omega^{\{k\}}$, and two $k$-subsets $U$ and $W$ are adjacent if $|U\cap W|=k-1$.
It is well-know that $\J(n,k)\cong \J(n,n-k)$ and $\J(n,1)\cong\K_n$.
Moreover, the Johnson graphs $\J(n,2)$ are also called the {\em Triangular graphs}.
Combining~\cite[Lemma 9.1.1]{BCN} and~\cite[Proposition 2.1]{JDLP}, we have the subsequent proposition.

\begin{proposition}\label{Johnson}
Let $\Ga=\J(n,k)$ be the Johnson graph with $2\leq k\leq [n/2]$. Then the following holds.
\begin{enumerate}[\rm (1)]
  \item $\Ga$ has order $\binom{n}{k}$, girth $3$, valency $k(n-k)$, diameter $d=\min\{k,n-k\}$,
        and the intersection array given by $b_i=(k-i)(n-k-i)$ and $c_i=i^2$ for $0\leq i\leq d$.
  \item $\Ga$ is geodesic transitive but not $2$-arc transitive, and $\Aut(\Ga)\cong\Sy_n\times \Sy_2$ if $n=2k$ and $\Aut(\Ga)\cong\Sy_n$ otherwise.
\end{enumerate}
\end{proposition}

Let $\Omega$ be a finite set with cardinality $2k+1$.
The {\em Odd graph $\O_k$} is the graph whose vertex set is the set of $k$-subsets of $\Omega$, and two $k$-subsets are adjacent if they are disjoint.
Then $\O_1\cong\C_3$, and for $k\geq2$, $\O_k$ has the following properties,
see~\cite[Propositions 9.1.7 and 9.1.8]{BCN},~\cite[Proposition 2.3]{JDLP} and \cite[Lemma 3.8]{JLX}.

\begin{proposition}\label{Odd}
Let $\Ga=\O_k$ be the Odd graph with $k\geq 2$. Then the following holds.
\begin{enumerate}[\rm (1)]
  \item $\Ga$ has order $\binom{2k+1}{k}$, valency $k+1$, diameter $k$,
        and the intersection array given by $b_i=k+1-[\frac{1}{2}(i+1)]$ and $c_i=[\frac{1}{2}(i+1)]$ for $0\leq i\leq d$ except that $b_k=0$.
  \item $\Ga$ is geodesic transitive and $3$-transitive, and $\Aut(\Ga)\cong\Sy_{2k+1}$.
  \item The double cover $2.\Ga$ of $\Ga$ is geodesic transitive and $3$-transitive, and $\Aut(2.\Ga)\cong\Sy_{2k+1}\times\ZZ_2$.
\end{enumerate}
\end{proposition}

Let $\F_q$ be a finite field of order $q$, and let $V$ be an $n$-dimensional vector space over $\F_q$.
Let $V_k$ be the set of $k$-subspaces of $q$.
Write $\genfrac{[}{]}{0pt}{1}{n}{k}$ the Gaussian $q$-binomial coefficient, that is,
\begin{equation*}
\genfrac{[}{]}{0pt}{1}{n}{k}=\frac{(q^n-1)(q^{n-1}-1)\cdots(q^{n-k+1-1}-1)}{(q^{m}-1)(q^{m-1}-1)\cdots(q-1)}.
\end{equation*}
The {\em Grassmann graph $\G_q^n(k)$} is the graph with vertex set $V_k$ such that two $k$-subspaces $U$ and $W$ are adjacent if $U\cap W$ has dimension $k-1$,
namely, $U\cap W$ is a hyperplane in $U$ and $W$.
Then $\G_q^n(k)\cong\G_q^n(n-k)$ (see~\cite[P. 268, Remarks]{BCN}) and $\G_q^n(1)\cong \G_q^n(n-1)\cong\K_{(q^n-1)/(q-1)}$.
Combing by~\cite[Theorems 9.3.1 and 9.3.3]{BCN} and~\cite[Corollary 2.9]{FH} we have the next result.

\begin{proposition}\label{Grassmann}
Let $\Ga=\G_q^d(k)$ be the Grassmann graph with $1<k<d-1$. Then the following holds.
\begin{enumerate}[\rm (1)]
  \item $\Ga$ has order $\genfrac{[}{]}{0pt}{1}{n}{k}$, diameter $d:=\min\{k,n-k\}$, and the intersection array given by
        \begin{equation*}
        b_i=q^{2i+1}\genfrac{[}{]}{0pt}{1}{k-i}{1}\genfrac{[}{]}{0pt}{1}{d-k-i}{1},
        c_i=\genfrac{[}{]}{0pt}{1}{i}{1}^2, \text{where~} 0\leq i\leq d.
        \end{equation*}
  \item $\Ga$ is geodesic transitive, $\Aut(\Ga)\cong\PGammaL(d,q)$ if $n\neq 2k$, and $\Aut(\Ga)\cong\PGammaL(d,q).\ZZ_2$ if $n=2k$.
\end{enumerate}
\end{proposition}

Let $n=2m+1$ be an odd integer and let $V$ be an $n$-dimensional vector space over $\F_q$.
The {\em doubled Grassmann graph} $2.\G_q^{2m+1}(m)$ has vertex set $V_m\cup V_{m+1}$,
where two vertices $U\in V_m$ and  $W\in V_{m+1}$ are adjacent if $U\subset W$.
This families of graphs are also studied in~\cite{FH}, and the authors showed that $2.\G_q^{2m+1}(m)$ is geodesic transitive of girth $6$,
and hence it is $3$-transitive.
Moreover, as noted in~\cite[P. 272]{BCN},
the doubled Grassmann graph $2.\G_q^{2m+1}(m)$ is the double cover of $\G_q^{2m+1}(m)$.
This leads to the following result.

\begin{proposition}\label{doubledGrassmann}
The doubled Grassmann graph $2.\G_q^{2m+1}(m)$ is geodesic transitive and $3$-transitive.
\end{proposition}

As defined in~\cite[P. 84]{GR}, a {\em generalized polygon}, or more precisely, a {\em generalized $d$-gon}, is a bipartite graph with diameter $d$ and girth $2d$.
Following~\cite{Weiss85}, we refer to the generalized polygons associated with the simple groups $\PSL_3(q)$, $\Sp_4(q)$ and $\G_2(q)$,
where $q$ is a prime power, as {\em classical generalized polygons}.
These are denoted by $\Delta_{3,q}$, $\Delta_{4,q}$ and $\Delta_{6,q}$ respectively.
The following result is quoted from~\cite[Section 6.5]{BCN} and~\cite[Lemma 2.5]{JDLP}.

\begin{proposition}\label{d-gon}
\begin{enumerate}[\rm (1)]
  \item $\Delta_{3,q}$ is a geodesic transitive and $4$-transitive graph with the intersection array $\{q+1,q,q;1,1,q+1\}$, and $\Aut(\Delta_{3,q})\cong\Aut(\PSL_3(q))$.
  \item $\Delta_{4,q}$ is a geodesic transitive and $5$-transitive graph with the intersection array $\{q+1,q,q,q;1,1,1,q+1\}$, and $\Aut(\Delta_{4,q})\cong\Aut(\Sp_4(q))$.
  \item $\Delta_{6,q}$ is a geodesic transitive and $7$-transitive graph with the intersection array $\{q+1,q,q,q,q,q;1,1,1,1,1,q+1\}$, and $\Aut(\Delta_{6,q})\cong\Aut(\G_2(q))$.
\end{enumerate}
\end{proposition}

Let $q=p^f$ be a prime power such that $q\equiv1\pmod{4}$, and let $\F_q$ be the finite field of order $q$.
The {\em Paley graph} $\P(q)$ is defined to be the graph with vertex set $\F_q$, and two vertices $u,v$ are adjacent if and only if $u-v$ is a nonzero square in $\F_q$.
This graph was first defined by Paley~\cite{Paley}, and $\P(q)=\Cay(\F_q^+,S)$ with $S=\la\lambda^2\ra$,
where $\F_q^+$ is the additive group of $\F_q$ and $\lambda$ be a primitive element of $\F_q$.
Furthermore,  $\P(q)$ has valency $(q-1)/2$ and diameter $2$, and $\P(q)$ has girth $3$ for $q>5$ and $\P(5)\cong\C_5$.
From~\cite[P. 221]{GR} and~\cite[Theorem 1.2]{JDLP}, we know that $\P(q)$ has the subsequent properties.

\begin{proposition}\label{Paley}
\begin{enumerate}[\rm (1)]
  \item  The Paley graph $\P(q)$ is distance transitive with intersection array $\{(q-1)/2,(q-1)/4;1,(q-1)/4\}$, and $\Aut(\Ga)\cong\ZZ_p^f:(\ZZ_{(q-1)/2}.\ZZ_f)$.
  \item  $\P(q)$ is geodesic transitive if and only if $q=5$ or $9$. Moreover, $\P(5)\cong \C_5$ and $\P(9)\cong H(2,3)$.
\end{enumerate}
\end{proposition}

Let $V$ be the affine space of $\AG(2,q)$.
Let $U=\{(x,y)\mid x,y\in\F_q\}$ be the set of points of $V$.
For each $m,k\in\F_q$, define
\begin{align*}
[m,k]=
\begin{cases}
\{(x,y): x,y\in\F_q\mid mx+y=k\}, \text{~if~} m\neq0,\\
\{(k,y): y\in\F_q\}, \text{~if~} m=0.
\end{cases}
\end{align*}
Then $[m,k]$ forms a translate of $1$-dimensional subspaces of $V$, which is an affine line of $V$.
Let $W=\{[m,k]:m,k\in\F_q\}$ be the set of affine lines of $V$.
Then each lines has $q$ points on it, and two distinct points of $U$ are joined by a unique line.
For given $m\in \F_q$, let $W_m:=\{[m,k]:k\in \F_q\}$.
Then $W_m$ is an partition of $V$, called a {\em parallel class}.
Now, the incidence structure $(U,W_m,W)$ is a resolvable divisible design $RGD(q,q,q,q^2)$ (see~\cite[P. 439]{BCN}).
Moreover, $(U,W_m,W)$ is also a $1$-$(q^2,q,q)$ design (see~\cite[P. 8]{JP}).
Define the graph $\mathcal{AG}(2,q)$ as the bipartite graph with bipartitions $U$ and $W$,
where a point $(x,y)$ is adjacent to a line $[m,k]$ if $(x,y)\in [m,k]$.
This graph was constructed as a distance transitive antipodal cover of the complete bipartite graph $\K_{q,q}$,
called {\em $\AG(2,q)$ minus a parallel class}, see~\cite{BCN}.
Clearly, $\mathcal{AG}(2,q)\cong\C_4$ and $\mathcal{AG}(2,3)$ is also known as the Pappus graph.

\begin{lemma}\label{AG2q}
Let $\Ga=\mathcal{AG}(2,q)$ with $q=p^f$ a prime power. Then $\Ga$ has the following properties.
\begin{enumerate}[\rm (1)]
  \item  $\Ga$ is a distance transitive graph with intersection array $\{q,q-1,q-1,1;1,1,q-1,q\}$.
  \item  $\Ga$ is an  antipodal cover of $\K_{q,q}$.
  \item  For each $a,c\in\F_q^*$ and $b,d\in\F_q$, define
         \begin{align*}
         &\tau_{a,b,c,d}:(x,y)\mapsto (ax+b,cy+d), \\
         &\a:(x,y) \mapsto [-x,-y] \text{~and~} [x,y]\mapsto (-x,-y), \text{~where~} (x,y)\in U.
         \end{align*}
         Then $\tau_{a,b,c,d}$ and $\a$ induce three automorphisms of $\Ga$, still denoted by $\tau_{a,b,c,d}$ and $\a$, respectively.
  \item  Let $G=\la\tau_{a,b,c,d},\a\mid a,c\in\F_q^*,b,d\in\F_q\ra$. Then $\Ga$ is $G$-geodesic transitive and $(G,3)$-transitive.
\end{enumerate}
\end{lemma}

\proof By~\cite[P. 273]{BCN}, we know that parts (1) and (2) of Lemma~\ref{AG2q} holds.

Recall that $U=\{(x,y)\mid x,y\in\F_q\}$ and $W=\{[m,k]:m,k\in\F_q\}$.
A direct computation may easily verify that
\begin{equation*}
[m,k]^{\tau_{a,b,c,d}}=[a^{-1}mc,a^{-1}mbc+ck+d].
\end{equation*}
Thus, $\tau_{a,b,c,d}$ and $\a$ are two permutations of $V(\Ga)=U\cup W$.
Moreover, for each edge $\{(x,y),[m,k]\}$ of $\Ga$, we have $(x,y)\in [m,k]$, and
\begin{align*}
&\{(x,y),[m,k]\}^{\tau_{a,b,c}}=\{(ax+b,cy+d),[a^{-1}mc,a^{-1}mbc+ck+d]\}\in E(\Ga),\\
&\{(x,y),[m,k]\}^{\a}=\{(-x,-y),[-m,-k]\}\in E(\Ga).
\end{align*}
It follows that $\tau_{a,b,c,d}$ and $\a$ induce two automorphisms of $\Ga$, as part (3).

To prove part (4), let $G=\la\tau_{a,b,c,d},\a\mid a,c\in\F_q^*,b,d\in\F_q\ra$.
Clearly, $\la\tau_{a,b,c,d}\mid a,c\in\F_q^*,b,d\in\F_q\ra$ acts transitively on both $U$ and $W$, and $\a$ interchanges $U$ and $W$.
Thus, $G$ is transitive on $V(\Ga)$.
Let
\begin{equation*}
u=[0,0], v=(0,0),w=[0,1]\text{~and~} z=(1,1).
\end{equation*}
Then $(u,v,w,z)$ is a $3$-geodesic of $\Ga$, and
\begin{equation*}
G_u=\la\tau_{a,b,c,0}\mid a,c\in\F_q^*,b\in\F_q\ra,~~G_{u,v}=\la \tau_{a,0,c,0}\mid a,c\in\F_q^*\ra, \text{~and~} G_{u,v,w}=\la\tau_{a,0,1,0}\mid a\in\F_q^*\ra.
\end{equation*}
Moreover, straightforward calculation shows that
\begin{align*}
&\Ga(u)=[0,0]=\{(0,i)\mid i \in\F_q\},~~
\Ga(v)\cap\Ga_2(u)=\{[0,i]\mid i\in \F_q^*\},\\
&\Ga(w)\cap\Ga_3(u)=[0,1]\setminus\{(0,0)\}=\{(i,1)\mid i\in\F_q^*\}.
\end{align*}
It follows that $G_u$, $G_{u,v}$ and $G_{u,v,w}$ act transitively on $\Ga(u)$, $\Ga(v)\cap\Ga_2(u)$ and $\Ga(w)\cap\Ga_3(u)$, respectively.
By Corollary~\ref{condition1} we know that $\Ga$ is $(G,3)$-geodesic transitive.
From part (1), $\Ga$ is distance transitive with intersection array $\{q,q-1,q-1,1;1,1,q-1,q\}$.
Thus, $\Ga$ has diameter $4$, girth $6$, and the parameter $b_4=1$.
It follows from Proposition~\ref{geo-tran} that $\Ga$ is $G$-geodesic transitive.
Since $\Ga$ has girth $6$, each $3$-geodesic is also a $3$-arc, and so $\Ga$ is $(G,3)$-arc transitive.
From~\cite[Proposition 17.2]{Biggs} we know that $\Ga$ is at most $3$-arc transitive.
Therefore, $\Ga$ is $(G,3)$-transitive, completing the proof. \qed

In~\cite[Section 7.5]{BCN}, the authors meticulously enumerate all distance transitive graphs with valency at most $13$.
Most of these graphs are already introduced as above
and the Webpage~\cite{BJWB}.
Consequently, we can expediently ascertain the geodesic transitivity of these graphs by using Magma~\cite{Magma}.
Next, we provide the detailed construction of graphs that are not included in the database~\cite{BJWB} or the graphs introduced above.

\begin{example}\label{graph1456-6}
Let $G=\Aut(\G_2(3))$ and let $H=3^2.(3\times 3_+^{1+2}):\D_8$ be a maximal subgroup of $G$.
Using Magma~\cite{Magma}, up to isomorphism, there exists a unique connected distance transitive graph of valency $6$, denoted by $\GG_{1456,6}$,
among the orbital graphs of $G$ acting on $[G:H]$.
Moreover, $\Aut(\GG_{1456,6})\cong G$ and $\GG_{1456,6}$ is geodesic transitive with intersection array  $\{6,3,3,3,3,3;1,1,1,1,1,2\}$.
\end{example}

\begin{example}\label{graph10-12}
Let $G=\Aut(^2\F_4(2)')$.
Let $H=2.[2^9]:5:4$ and let $K=2^2.[2^9].\Sy_3$ be two maximal subgroups of $G$.
By Magma~\cite{Magma}, up to isomorphism, there exist the unique connected distance transitive graph of valency $10$ and $12$,
denoted by $\GG_{1755,10}$ and $\GG_{2925,12}$, respectively,
among the orbital graphs of $G$ acting on $[G:H]$ and $[G:K]$.
Moreover, the following properties holds:
\begin{enumerate}[\rm (1)]
  \item $\GG_{1755,10}$ is a geodesic transitive graph with intersection array $\{10,8,8,8;1,1,1,5\}$ and $\Aut(\GG_{1755,10})\cong G$;
  \item $\GG_{2925,12}$ is a geodesic transitive graph with intersection array  $\{12,8,8,8;1,1,1,3\}$ and $\Aut(\GG_{2925,12})\cong G$.
\end{enumerate}
\end{example}

The adjacency matrices of the graph $\mathcal{AG}(2,q)$ for $q\leq 9$ are provided in \cite{BJWB}.
To facilitate calculations with Magma~\cite{Magma}, we present the structure of $\mathcal{AG}(2,q)$ for $q=11$ or $13$ using Cayley graphs.

\begin{example}\label{AG211}
Let $G=\la a,b,c\mid a^{11}=b^{11}=c^2=[a,b]=1,a^c=a^{-1},b^c=b^{-1}\ra$ be the group of order $2\cdot 11^2$.
Let $\Ga=\Cay(G,S)$ with
\begin{equation*}
S=\{a^{10}b^{10}c,a^9b^{10}c,a^2b^2c,a^4bc,b^5c,ab^6c,a^3b^4c,a^6b^2c,a^8b^5c,a^5b^4c,a^7b^6c\}.
\end{equation*}
By Magma~\cite{Magma}, $\Ga\cong\mathcal{AG}(2,11)$ is a geodesic transitive and $3$-transitive graph.
Moreover, $\Aut(\mathcal{AG}(2,11))\cong 11^{1+3}_+:(\ZZ_5^2.\D_8)$.
\end{example}

\begin{example}\label{AG213}
Let $G=\la a,b,c\mid a^{13}=b^{13}=c^2=[a,b]=1,a^c=a^{-1},b^c=b^{-1}\ra$ be the group of order $2\cdot 13^2$.
Let $\Ga=\Cay(G,S)$ with
\begin{equation*}
S=\{b^{12}c,a^{12}b^{4}c,a^9bc,a^4bc,ab^4c,a^{10}b^5c,a^{11}b^6c,a^3b^5c,a^2b^6c,a^7b^{10}c,a^6b^{10}c,a^8b^7c,a^5b^7c\}.
\end{equation*}
By Magma~\cite{Magma}, $\Ga\cong\mathcal{AG}(2,13)$ is a geodesic transitive and $3$-transitive graph.
Moreover, $\Aut(\mathcal{AG}(2,13))\cong 13^{1+3}_+:(3^2.2^{2+3})$.
\end{example}

In our search, we have identified exactly four graphs, in addition to the Paley graph,
which are distance transitive but not geodesic transitive.
For convenience, we also provide the construction of these graphs.

\begin{example}\label{graph22-6}
Let $G=\la a,b\mid a^{11}=b^2=1, a^b=a^{-1}\ra\cong\D_{22}$ and let $S=\{b,a^2b,a^6b,a^{7}b$, $a^{8}b, a^{10}b\}$.
Define $\GG_{22,6}=\Cay(G,S)$.
Then $\GG_{22,6}$ is distance transitive with intersection array $\{6,5,3;1,3,6\}$, and $\Aut(\GG_{22,6})\cong\PGL_2(11)$.
\end{example}

We remark that $\GG_{22,6}$ is in fact the incidence graph of $2$-$(11,6,3)$ design.

\begin{example}\label{graph64-8}
Let $G=\la a,b,c,d\ra$ denote a nonabelian group of order $64$ with the following relations:
\begin{align*}
a^2=b^4=c^4=d^2=1,~b^a=b^{-1},~c^a=c^{-1},~[a,d]=[b,c]=[b,d]=[c,d]=1.
\end{align*}
Let $S=\{a,ab,abc,abd,ab^3c^3d,ab^2c^3,ac^2d,ac^3d\}$ and let $\GG_{64,8}:=\Cay(G,S)$.
Then $\GG_{64,8}$ is distance transitive with intersection array $\{8, 7, 6, 1; 1, 2, 7, 8\}$, and $\Aut(\GG_{64,8})=2^{2+6}:\AGL_1(7)$.
\end{example}

\begin{example}\label{graph280-9}
Let $T=\PSL_3(4)$ and let $A=\Aut(T)\cong T.(2\times\Sy_3)$.
Then $A$ has a maximal subgroup $H\cong 3^2:2\Sy_4\times 2$ with index $280$.
By Magma~\cite{Magma}, up to isomorphism, there exists a unique connected distance transitive graph of valency $9$, denoted by $\GG_{280,9}$,
among the orbital graphs of $A$ acting on $[A:H]$.
Moreover, $\GG_{280,9}$ has intersection array  $\{9,8,6,3;1,1,3,8\}$ and $\Aut(\GG_{280,9})\cong A$.
\end{example}

\begin{example}\label{graph68-12}
Let $G=\PGammaL_2(16)$.
Then $G$ has a subgroup $H\cong\A_6.\ZZ_2^2$ with index $68$.
Using Magma~\cite{Magma}, it has been established that, up to isomorphism, there exist the unique connected distance transitive graph with valency $12$,
denoted as $\GG_{68,12}$, among the orbital graphs of $G$ acting on $[G:H]$.
Furthermore, $\GG_{68,12}$ has intersection array $\{12,10,3;1,3,8\}$ and $\Aut(\GG_{68,12})\cong G$.
\end{example}

\section{Geodesic transitive graphs of small valency}\label{sum}

Using Lemma~\ref{condition}, we can employ Magma~\cite{Magma} to determine the $s$-geodesic transitivity of all distance transitive graphs with valency at most $13$.
The Magma calculator program is provided in subsection~\ref{magmacode}.
Through computation, we obtain the following results:

\begin{lemma}\label{ngeo-trans}
\begin{enumerate}[\rm (1)]
  \item  The Paley graph $\P(q)$, with $q=13$, $17$, or $25$, is $1$-geodesic transitive but not $2$-geodesic transitive.
  \item The graphs $\GG_{22,6}$, $\GG_{64,8}$ and $\GG_{68,12}$ are $2$-geodesic transitive but not $3$-geodesic transitive.
  \item The graph $\GG_{280,9}$ is $3$-geodesic transitive but not $4$-geodesic transitive.
\end{enumerate}
\end{lemma}

All calculation results are presented in Tables~\ref{val3}--\ref{val13}, leading to the verification of Theorem~\ref{geo-small}.
Noting that every geodesic transitive graph of valency $1$ or $2$ is a cycle $\C_n$ with $n\geq2$.
Using this fact and combining it with Theorem~\ref{geo-small}, we provide the classification of all geodesic transitive graphs of valency at most $13$.

\subsection{MAGMA Code}\label{magmacode}~

For a graph $X$, the \texttt{IsDistanceTransitive(X)} and \texttt{IntersectionArray(X)} procedures in Magma~\cite{Magma}
can be used to determine whether the graph is distance transitive and to compute its intersection array, respectively.
In the following Magma~\cite{Magma} code, it will determine whether the graph $X$ is $s$-geodesic transitive.

\medskip
\noindent{}{\bf Input:} a finite graph $X$ \\
{\bf Output:} the graph $X$ is geodesic transitive or it is $s$-geodesic transitive but not $(s+1)$-geodesic transitive.
\begin{verbatim}
geodesic:=function(X);
f:=IsVertexTransitive(X);
if f eq true then
   d:=Diameter(X);
   A:=AutomorphismGroup(X);
   V:=VertexSet(X);
   V1:=Vertices(X);
   u:=V.1;
   S:=Sphere(u,d);
   Geo:=[Geodesic(u, i): i in S];
   geo:=Geo[1];
   x:=[Index(V1,geo[i]): i in [1..#geo]];
   F:=[];
   for i in [1..d] do
        sgeo:=[x[j]: j in [1..i]];
        Ssgeo:=Stabilizer(A,sgeo);
        us:=x[i];
        N:=Neighbours(V!us) meet Sphere(u,i);
         a:=Index(V1,SetToSequence(N)[1]);
         D:={V!i: i in a^Ssgeo};
         f1:=N eq D;
         if f1 eq true then
            Append(~F,1);
         else
            Append(~F,0);
         end if;
   end for;
   F1:=[i: i in [1..#F] | F[i] eq 0];
   if #F1 eq 0 then
       return"the graph X is geodesic transitive";
   else
     return"the graph X is", F1[1]-1, "-geodesic transitive but not",
     F1[1], "-geodesic transitive";
   end if;
else
   return"the graph X is not vertex transitive";
end if;
end function;
\end{verbatim}

\subsection{The Tables}
We are now ready to display the list of all connected geodesic transitive graphs of valency $3$ to $13$.
The notation for most graphs listed in the tables of this paper is introduced in Section~\ref{example}.
For the reminder of the graphs, we use $\GG_{n,k}$ to denote a distance transitive graph of order $n$ and valency $k$,
and use $\GG_{n,k}^1$, $\GG_{n,k}^2,\ldots$ to denote non-isomorphic graphs.
For a graph $\Ga$, denoted by $d$ and $g$ the diameter and girth of $\Ga$,
and by $A$ and $A_u$ the full automorphism group of $\Ga$ and the stabilizer $A_u$ of a vertex $u$ in $A$, respectively.
Moreover, $s$ indicates that the graph is $s$-arc transitive but not $(s+1)$-arc transitive.
The column ``GT" provides the information on whether the graph is geodesic transitive.
Furthermore, for the notation of groups in the Tables~\ref{val3}--\ref{val13}, we refer to Atlas~\cite{Atlas}.

\newpage
{\fontsize{9pt}{11pt}\selectfont
\begin{table}[!htbp]
\caption{The distance transitive graphs of valency $3$}
\label{val3}
\begin{tabular}{l|l|c|c|c|c} \hline
$\Ga$ & Intersection array & $(d,g,s)$ & $A$ & $A_u$ & GT  \\ \hline
  $\K_4$   & $\{3; 1\}$    & $(1,3,2)$ & $\Sy_4$ & $\Sy_3$ & Y    \\
  $\K_{3,3}$  & $\{3,2;1,3\}$   & $(2,4,3)$ & $\Sy_3\wr\Sy_2$ & $\Sy_3\times\Sy_2$& Y   \\
  $\GG_{10,3}$ & $\{3,2;1,1\}$   & $(2,5,3)$  & $\Sy_5$ & $\Sy_3\times\Sy_2$ &Y   \\
  $H(3,2)$ & $\{3,2,1; 1,2,3\}$ & $(3,4,2)$ & $\ZZ_2^3\rtimes\Sy_3$ & $\Sy_3$ & Y   \\
  $\GG_{14,3}$ & $\{3,2,2; 1,1,3\}$ & $(3,6,4)$ & $\PGL_2(7)$ & $\Sy_4$ & Y   \\
  $\mathcal{AG}(2,3)$ & $\{3,2,2,1; 1,1,2,3\}$ & $(4,6,3)$ & $\ZZ_3.(\Sy_3\wr\Sy_2)$& $\Sy_3\times\Sy_2$ & Y  \\
  $\GG_{28,3}$ & $\{3,2,2,1;1,1,1,2\}$ & $(4,7,3)$ & $\PGL_2(7)$ & $\Sy_3\times\Sy_2$ & Y \\
  $\GG_{30,3}$ & $\{3,2,2,2;1,1,1,3\}$ & $(4,8,5)$ & $\PGammaL_2(9)$ & $\Sy_4\times\Sy_2$ & Y   \\
  $\GG_{20,3}^1$ & $\{3,2,1,1,1;1,1,1,2,3\}$ & $(5,5,2)$ & $\A_5\times\ZZ_2$ & $\Sy_3$ & Y    \\
  $\GG_{20,3}^2$ & $\{3,2,2,1,1;1,1,2,2,3\}$ & $(5,6,3)$ & $\Sy_5\times\ZZ_2$ & $\Sy_3\times\Sy_2$ & Y \\
  $\GG_{102,3}$ & $\{3,2,2,2,1,1,1; 1,1,1,1,1,1,3\}$ & $(6,9,4)$ & $\mathrm{L}_2(17)$ & $\Sy_4$ & Y  \\
  $\GG_{90,3}$ & $\{3,2,2,2,2,1, 1,1;1,1,1,1,2,2,2,3\}$ & $(8,10,5)$ & $\ZZ_3.\PGammaL_2(9)$ & $\Sy_4\times\Sy_2$& Y  \\ \hline
\end{tabular}
\end{table}

\begin{table}[!htbp]
\centering
\caption{The distance transitive graphs of valency $4$}
\label{val4}
\begin{tabular}{l|l|c|c|c|c} \hline
  $\Ga$ & Intersection array & $(d,g,s)$ & $A$ & $A_u$ & GT  \\ \hline
  $\K_4$ & $\{4;1\}$ & $(1,3,2)$ & $\Sy_5$ & $\Sy_4$ & Y   \\
  $\K_{4,4}$ & $\{4,3;1,4\}$  & $(2,4,3)$ & $\Sy_4\wr\Sy_2$ & $\Sy_3\times\Sy_4$ & Y \\
  $\K_{5,5}-5\K_2$ & $\{4,3,1;1,3,4\}$  & $(3,4,2)$ & $\Sy_5\times\Sy_2$ & $\Sy_4$ & Y   \\
  $H(4,2)$ & $\{4,3,2,1;1,2,3,4\}$  & $(4,4,2)$ & $\ZZ_2^4:\Sy_4$ & $\Sy_4$ & Y    \\
  $\O_4$ & $\{4,3,3;1,1,2\}$  & $(3,6,3)$ & $\Sy_7$ & $\Sy_4\times\Sy_3$ & Y    \\
  $2.\O_4$ & $\{4,3,3,2,2,1,1;1,1,2,2,3,3,4\} $ & $(7,6,3)$ & $\Sy_7\times\ZZ_2$ & $\Sy_4\times\Sy_3$ & Y    \\
  $\K_{3[2]}$ & $\{4,1;1,4\}$ & $(2,3,1)$ & $\Sy_2\wr\Sy_3$ & $\D_8$& Y   \\
  $H(2,3)$ & $\{4,2;1,2\}$ & $(2,3,1)$ & $\Sy_3\wr\Sy_2$ & $\D_8$& Y   \\
  $\GG_{15,4}$ & $\{4,2,1;1,1,4\}$ & $(3,3,1)$ & $\Sy_5$ & $\D_8$& Y   \\
  $\GG_{21,4}$ & $\{4,2,2;1,1,2\}$ & $(3,3,1)$ & $\PGL_2(7)$ & $\D_{16}$& Y  \\
  $\GG_{14,4}$ & $\{4,3,2;1,2,4\}$ & $(3,4,2)$ & $\PGL_2(7)$ & $\Sy_4$ & Y   \\
  $\Delta_{3,3}$ & $\{4,3,3;1,1,4\}$ & $(3,6,4)$ & $\mathrm{L}_3(3).\ZZ_2$ & $3^2.\Q_8.\Sy_3$ & Y  \\
  $\GG_{45,4}$ & $\{4,2,2,2;1,1,1,2\}$ & $(4,3,1)$ & $\PGammaL_2(9)$ & $[2^5]$ & Y  \\
  $\mathcal{AG}(2,4)$ & $\{4,3,3,1;1,1,3,4\}$ & $(4,6,3)$ & $\ZZ_2^2.((\Sy_4\times\A_4).\ZZ_2)$ & $(3\times\A_4).2$ & Y   \\
  $\Delta_{5,3}$ & $\{4,3,3,3,3,3;1,1,1,1,1,4\}$ & $(6,12,7)$ & $\G_2(3).\ZZ_2$ & $[3^5]:\GL_2(3)$ & Y  \\ \hline
\end{tabular}
\end{table}

\begin{table}[!htbp]
\centering
\caption{The distance transitive graphs of valency $5$}
\label{val5}
\begin{tabular}{l|l|c|c|c|c} \hline
  $\Ga$ & Intersection array & $(d,g,s)$ & $A$ & $A_u$ & GT  \\ \hline
  $\K_6$ & $\{5;1\}$ & $(1,3,2)$ & $\Sy_6$ & $\Sy_5$ & Y   \\
  $\K_{5,5}$ & $\{5,4;1,5\}$  & $(2,4,3)$ & $\Sy_5\wr\Sy_2$ & $\Sy_4\times\Sy_5$ & Y  \\
  $\K_{6,6}-6\K_2$ & $\{5,4,1;1,4,5\}$  & $(3,4,2)$ & $\Sy_6\times\Sy_2$ & $\Sy_5$ & Y   \\
  $H(5,2)$ & $\{5,4,3,2,1;1,2,3,4,5\}$  & $(5,4,2)$ & $\ZZ_2^5:\Sy_5$ & $\Sy_5$ & Y   \\
  $\square_5$ & $\{5,4;1,2\}$  & $(2,4,2)$ & $\ZZ_2^4:\Sy_5$ & $\Sy_5$ & Y    \\
  $\O_5$ & $\{5,4,4,3;1,1,2,2\}$  & $(4,6,3)$ & $\Sy_9$ & $\Sy_5\times\Sy_4$ & Y   \\ \hline
  $2.\O_5$ & $\{5,4,4,3,3,2,2,1,1; $ & $(9,6,3)$ & $\Sy_9\times\ZZ_2$ & $\Sy_5\times\Sy_4$ & Y   \\
           & \hspace{0.5em}$1,1,2,2,3,3,4,4,5\}$  & & & & \\\hline
  $\GG_{12,5}$ & $\{5,2,1;1,2,5\}$ & $(3,3,1)$ & $\A_5\times\ZZ_2$ & $\D_{10}$ & Y \\
  $\GG_{36,5}$ & $\{5,4,2;1,1,4\}$ & $(3,5,2)$ & $\Sy_6.\ZZ_2$ & $\AGL_1(5)\times\ZZ_2$ & Y  \\
  $\GG_{22,5}$ & $\{5,4,3;1,2,5\}$ & $(3,4,2)$ & $\PGL_2(11)$ & $\A_5$  & Y  \\
  $\Delta_{3,4}$ & $\{5,4,4;1,1,5\}$ & $(3,6,4)$ & $\mathrm{L}_3(4).(\ZZ_2\times\Sy_3)$& $[4^2]:\GL_2(4)$ & Y  \\
  $\GG_{32,5}$ & $\{5,4,1,1;1,1,4,5\}$ & $(4,5,2)$ & $2_{-}^{1+4}:\A_5$ & $\A_5$ & Y  \\
  $\mathcal{AG}(2,5)$ & $\{5,4,4,1;1,1,4,5\}$ & $(4,6,3)$ & $\ZZ_5.((\F_{20}\times\F_{20}):\ZZ_2)$ & $\AGL_1(5)\times\ZZ_4$ & Y   \\
  $\Delta_{4,4}$ & $\{5,4,4,4;1,1,1,5\}$ & $(4,8,5)$ & $\PSp_4(4).\ZZ_4$& $[4^3]:\GammaL_2(4)$ & Y \\  \hline
\end{tabular}
\end{table}

\newpage
\begin{table}[!htbp]
\centering
\caption{The distance transitive graphs of valency $6$}
\label{val6}
\begin{tabular}{l|l|c|c|c|c} \hline
  $\Ga$ & Intersection array & $(d,g,s)$ & $A$ & $A_u$ & GT\\ \hline
  $\K_7$  & $\{6;1\}$ & $(1,3,2)$ & $\Sy_7$ & $\Sy_6$ & Y  \\
  $\K_{6,6}$ & $\{6,5;1,6\}$  & $(2,4,3)$ & $\Sy_6\wr\Sy_2$ & $\Sy_5\times\Sy_6$ & Y   \\
  $\K_{7,7}-7\K_2$ & $\{6,5,1;1,5,6\}$  & $(3,4,2)$ & $\Sy_7\times\Sy_2$ & $\Sy_6$ & Y   \\
  $H(6,2)$ & $\{6,5,4,3,2,1;1,2,3,4,5,6\}$ & $(6,4,2)$ & $\ZZ_2^6:\Sy_6$ & $\Sy_6$ & Y    \\
  $\square_6$ & $\{6,5,4;1,2,6\}$  & $(3,4,2)$ & $\ZZ_2^5:\Sy_6$ & $\Sy_6$ & Y  \\
  $\O_6$ & $\{6,5,5,4,4;1,1,2,2,3\}$  & $(5,6,3)$ & $\Sy_{11}$ & $\Sy_6\times\Sy_5$ & Y    \\ \hline
  $2.\O_6$ & $\{6,5,5,4,4,3,3,2,2,1,1; $ & $(11,6,3)$ & $\Sy_{11}\times\ZZ_2$ & $\Sy_6\times\Sy_5$ & Y    \\
           & \hspace{0.5em}$1,1,2,2,3,3,4,4,5,5,6\}$  & & &  \\\hline
  $\K_{4[2]}$  & $\{6,1;1,6\}$ & $(2,3,1)$ & $\Sy_2\wr\Sy_4$ & $\ZZ_2\wr\Sy_3$ & Y  \\
  $\K_{3[3]}$ & $\{6,2;1,6\}$ & $(2,3,1)$ & $\Sy_3\wr\Sy_3$ & $\ZZ_2\times(\Sy_3\wr\Sy_2)$ & Y  \\
  $\J(5,2)$ & $\{6,2;1,4\}$ & $(2,3,1)$ & $\Sy_5$ & $\Sy_3\times\ZZ_2$  & Y  \\
  \rowcolor{yellow}
  $\P(13)$ & $\{6,3;1,3\}$ & $(2,3,1)$ & $\ZZ_{13}:\ZZ_6$ & $\ZZ_6$ & N   \\
  $\GG_{15,6}$ & $\{6,4;1,3\}$ & $(2,3,1)$ & $\Sy_6$ & $\Sy_4\times\ZZ_2$ & Y \\
  $H(2,4)$ & $\{6,3;1,2\}$ & $(2,3,1)$ & $\Sy_4\wr\Sy_2$ & $\Sy_3\wr\Sy_2$ & Y  \\
  $\GG_{52,6}$ & $\{6,3,3;1,1,2\}$ & $(3,3,1)$ & $\mathrm{L}_3(3).\ZZ_2$ & $3_{+}^{1+2}:\D_8$ & Y  \\
  $H(3,3)$ & $\{6,4,2;1,2,3\}$ & $(3,3,1)$ & $\Sy_3\wr\Sy_3$ & $\Sy_2\wr\Sy_3$ & Y  \\\hline
  $\GG_{63,6}^1$ & $\{6,4,4;1,1,3\}$ & $(3,3,1)$ & $\mathrm{U}_3(3).\ZZ_2$ & $4^2:\D_{12}$  & Y \\
  $\GG_{63,6}^2$ & $\{6,4,4;1,1,3\}$ & $(3,3,1)$ & $\mathrm{U}_3(3).\ZZ_2$ & $4\cdot\Sy_4:2$ & Y    \\ \hline
  $\GG_{42,6}$ & $\{6,5,1;1,1,6\}$ & $(3,5,2)$ & $\Sy_7$ & $\Sy_5$ & Y  \\
  $\GG_{57,6}$ & $\{6,5,2; 1,1,3\}$ & $(3,5,2)$ & $\mathrm{L}_2(19)$ & $\A_5$ & Y  \\
  \rowcolor{yellow}
  $\GG_{22,6}$ & $\{6,5,3; 1,3,6\}$ & $(3,4,2)$ & $\PGL_2(11)$ & $\A_5$  & N  \\
  $\Delta_{3,5}$ & $\{6,5,5;1,1,6\}$ & $(3,6,4)$ & $\mathrm{L}_3(5).\ZZ_2$ & $5^2:\GL_2(5)$ & Y  \\
  $\GG_{45,6}$ & $\{6,4,2,1; 1,1,4,6\}$ & $(4,3,1)$ & $\ZZ_3.\Sy_6$ & $\Sy_4\times\ZZ_2$  & Y  \\
  $\GG_{36,6}$ & $\{6,5,4,1; 1,2,5,6\}$ & $(4,4,2)$ & $\ZZ_3.\Sy_6.\ZZ_2$ & $\Sy_5$ & Y  \\
  $\GG_{1456,6}$ & $\{6,3,3,3,3,3; 1,1,1,1,1,2\}$ & $(6,3,1)$ & $\G_2(3).\ZZ_2$ & $3^2.(3\times3_{+}^{1+2}):\D_8$ & Y \\  \hline
\end{tabular}
\end{table}

\begin{table}[!htbp]
\centering
\caption{The distance transitive graphs of valency $7$}
\label{val7}
\begin{tabular}{l|l|c|c|c|c} \hline
  $\Ga$ & Intersection array & $(d,g,s)$ & $A$ & $A_u$ & GT  \\ \hline
  $\K_8$ & $\{7;1\}$ & $(1,3,2)$ & $\Sy_8$ & $\Sy_7$ & Y  \\
  $\K_{7,7}$ & $\{7,6;1,7\}$  & $(2,4,3)$ & $\Sy_7\wr\Sy_2$ & $\Sy_6\times\Sy_7$ & Y   \\
  $\K_{8,8}-8\K_2$ & $\{7,6,1;1,6,7\}$  & $(3,4,2)$ & $\Sy_8\times\Sy_2$ & $\Sy_8$ & Y   \\\hline
  $H(7,2)$ & $\{7,6,5,4,3,2,1;$ & $(7,4,2)$ & $\ZZ_2^7:\Sy_7$ & $\Sy_7$ & Y    \\
  & \hspace{0.5em}$1,2,3,4,5,6,7\}$  & & &  \\\hline
  $\square_7$ & $\{7,6,5;1,2,3\}$  & $(3,4,2)$ & $\ZZ_2^6:\Sy_7$ & $\Sy_7$ & Y   \\ \hline
  $\O_7$ & $\{7,6,6,5,5,4;$  & $(6,6,3)$ & $\Sy_{13}$ & $\Sy_7\times\Sy_6$ & Y   \\
                 & \hspace{0.5em}$1,1,2,2,3,3\}$  & & &  \\ \hline
  $2.\O_7$ & $\{7,6,6,5,5,\cdots,1,1; $ & $(13,6,3)$ & $\Sy_{13}\times\ZZ_2$ & $\Sy_7\times\Sy_6$ & Y   \\
           & \hspace{0.5em}$1,1,2,2,\cdots,6,6,7\}$  & & &  \\\hline
  $\GG_{50,7}$ & $\{7,6;1,1\}$ & $(2,5,3)$ & $\mathrm{U}_3(5).\ZZ_2$ & $\Sy_7$ & Y  \\
  $\GG_{30,7}$ & $\{7,6,4;1,3,7\}$ & $(3,4,2)$ & $\mathrm{L}_4(2).\ZZ_2$ & $2^3:\mathrm{L}_3(2)$ & Y  \\
  $\GG_{330,7}$ & $\{7,6,4,4;1,1,1,6\}$ & $(4,5,2)$ & $\M_{22}.\ZZ_2$ & $2^3:\mathrm{L}_3(2)\times\ZZ_2$  & Y \\
  $\mathcal{AG}(2,7)$ & $\{7,6,6,1;1,1,6,7\}$ & $(4,6,3)$ & $7^{1+2}:(\ZZ_6.\D_{12})$ & $\F_{42}\times\ZZ_6$  & Y  \\
  $\GG_{100,7}$ & $\{7,6,6,1,1;1,1,6,6,7\}$ & $(5,6,3)$ & $\PSigmaU_3(5)\times\ZZ_2$ & $\Sy_7$ & Y  \\
  $2.\G_2^5(2)$ & $\{7,6,6,4,4; 1,1,3,3,7\}$ & $(5,6,3)$ & $\mathrm{L}_5(2).\ZZ_2$ & $\ZZ_2^6:(\Sy_3\times\mathrm{L}_3(2))$ & Y \\ \hline
  $\GG_{990,7}$ & $\{7,6,4,4,4,1,1,1;$ & $(8,5,2)$ & $3.\M_{22}:2$ & $\ASL_3(2)\times\ZZ_2$ & Y   \\
                & \hspace{0.5em}$1,1,1,2,4,4,6,7\}$  & & & & \\ \hline
\end{tabular}
\end{table}

\newpage
\begin{table}[!htbp]
\centering
\caption{The distance transitive graphs of valency $8$}
\label{val8}
\begin{tabular}{l|l|c|c|c|c} \hline
  $\Ga$ & Intersection array & $(d,g,s)$ & $A$ & $A_u$ & GT \\ \hline
  $\K_9$ & $\{8;1\}$ & $(1,3,2)$ & $\Sy_9$ & $\Sy_8$ & Y  \\
  $\K_{8,8}$ & $\{8,7;1,8\}$  & $(2,4,3)$ & $\Sy_8\wr\Sy_2$ & $\Sy_8\times\Sy_7$ & Y  \\
  $\K_{9,9}-9\K_2$& $\{8,7,1;1,7,8\}$  & $(3,4,2)$ & $\Sy_9\times\Sy_2$ & $\Sy_8$ & Y  \\\hline
  $H(8,2)$  & $\{8,7,6,5,4,3,2,1;$ & $(8,4,2)$ & $\ZZ_2^8:\Sy_8$ & $\Sy_8$ & Y   \\
  & \hspace{0.5em}$1,2,3,4,5,6,7,8\}$  & & & \\\hline
  $\square_8$ & $\{8,7,6,5;1,2,3,8\}$  & $(4,4,2)$ & $\ZZ_2^7:\Sy_8$ & $\Sy_8$ & Y   \\ \hline
  $\O_8$ & $\{8,7,7,6,6,5,5;$  & $(7,6,3)$ & $\Sy_{15}$ & $\Sy_8\times\Sy_7$ & Y    \\
   & \hspace{0.5em}$1,1,2,2,3,3,4\}$  & & &  \\\hline
  $2.\O_8$  & $\{8,7,7,6,6,\cdots,1,1; $ & $(15,6,3)$ & $\Sy_{15}\times\ZZ_2$ & $\Sy_8\times\Sy_7$ & Y  \\
           & \hspace{0.5em}$1,1,2,2,\cdots,7,7,8\}$  & & &  \\\hline
  $\K_{5[2]}$ & $\{8,1;1,8\}$ & $(2,3,1)$ & $\Sy_2\wr\Sy_5$ & $\Sy_2\wr\Sy_4$ & Y  \\
  $\K_{3[4]}$ & $\{8,3;1,8\}$ & $(2,3,1)$ & $\Sy_4\wr\Sy_3$ & $\Sy_3\times(\Sy_4\wr\Sy_2)$ & Y  \\
  $\J(6,2)$ & $\{8,3;1,4\}$ & $(2,3,1)$ & $\Sy_6$ & $\Sy_4\times\ZZ_2$ & Y \\
  \rowcolor{yellow}
  $\P(17)$ & $\{8,4;1,4\}$ & $(2,3,1)$ & $\ZZ_{17}:\ZZ_{8}$ & $\ZZ_8$ & N   \\
  $H(2,5)$ & $\{8,4;1,2\}$ & $(2,3,1)$ & $\Sy_5\wr\Sy_2$& $\Sy_4\wr\Sy_2$ & Y  \\
  $\GG_{105,8}$ & $\{8,4,4;1,1,2\}$ & $(3,3,1)$ & $\mathrm{L}_3(4).\D_{12}$ & $2^{2+4}.3^2.2^2$ & Y  \\
  $\GG_{27,8}$ & $\{8,6,1;1,3,8\}$ & $(3,3,1)$ & $3^{1+2}_{+}:\GL_2(3)$& $\GL_2(3)$ & Y  \\
  $\GG_{30,8}$ & $\{8,7,4;1,4,8\}$ & $(3,4,2)$ & $\mathrm{L}_4(2).\ZZ_2$ & $2^3:\mathrm{L}_3(2)$ & Y  \\
  $\Delta_{3,7}$ & $\{8,7,7;1,1,8\}$ & $(3,6,4)$ & $\mathrm{L}_3(7).\Sy_3$ & $7^2:\GL_2(7)$  & Y   \\
  $\GG_{425,8}$ & $\{8,4,4,4;1,1,1,2\}$ & $(4,3,1)$ &  $\PSp_4(4).\ZZ_4$& $(2^2\times2^{2+4}:3):12$  & Y   \\
  $H(4,3)$  & $\{8,6,4,2;1,2,3,4\}$ & $(4,3,1)$ & $\Sy_3\wr\Sy_4$ & $\Sy_2\wr\Sy_4$ & Y  \\
  $\GG_{32,8}$ & $\{8,7,4,1;1,4,7,8\}$ & $(4,4,2)$ & $2^{1+6}:\PGL_3(2)$ & $\AGL_3(2)$ & Y \\
  \rowcolor{yellow}
  $\GG_{64,8}$ & $\{8,7,6,1;1,2,7,8\}$ & $(4,4,2)$ & $2^{2+6}:\AGL_1(7)$ & $\ZZ_2^3:(\ZZ_7:\ZZ_3)$ & N   \\
  $\mathcal{AG}(2,8)$ & $\{8,7,7,1;1,1,7,8\}$ & $(4,6,3)$ & $2^{3+6}:(\ZZ_7^2.\ZZ_6)$ & $\ZZ_2^3:(\ZZ_7^2.\ZZ_3)$  & Y \\ \hline
\end{tabular}
\end{table}

\begin{table}[!htbp]
\centering
\caption{The distance transitive graphs of valency $9$}
\label{val9}
\begin{tabular}{l|l|c|c|c|c} \hline
  $\Ga$ & Intersection array & $(d,g,s)$ & $A$ & $A_u$ & GT \\ \hline
  $\K_{10}$ & $\{9;1\}$ & $(1,3,2)$ & $\Sy_{10}$ & $\Sy_9$ & Y   \\
  $\K_{9,9}$ & $\{9,8;1,9\}$  & $(2,4,3)$ & $\Sy_9\wr\Sy_2$ & $\Sy_9\times\Sy_8$ & Y \\
  $\K_{10,10}-10\K_2$ & $\{9,8,1;1,8,9\}$  & $(3,4,2)$ & $\Sy_{10}\times\Sy_2$ & $\Sy_9$ & Y   \\\hline
  $H(9,2)$ & $\{9,8,7,6,5,4,3,2,1;$ & $(9,4,2)$ & $\ZZ_2^9:\Sy_9$ & $\Sy_9$ & Y   \\
  & \hspace{0.5em}$1,2,3,4,5,6,7,8,9\}$  & & &  \\\hline
  $\square_9$ & $\{9,8,7,6;1,2,3,4\}$  & $(4,4,2)$ & $\ZZ_2^8:\Sy_9$ & $\Sy_9$ & Y    \\ \hline
  $\O_9$ & $\{9,8,8,7,7,6,6,5;$  & $(8,6,3)$ & $\Sy_{17}$ & $\Sy_9\times\Sy_8$ & Y   \\
  & \hspace{0.5em}$1,1,2,2,3,3,4,4\}$  & & & & \\\hline
  $2.\O_9$ & $\{9,8,8,7,7,\cdots,1,1; $ & $(17,6,3)$ & $\Sy_{17}\times\ZZ_2$ & $\Sy_9\times\Sy_8$ & Y    \\
           & \hspace{0.5em}$1,1,2,2,\cdots,8,8,9\}$  & & & &  \\\hline
  $\K_{4[3]}$ & $\{9,2;1,9\}$ & $(2,3,1)$ & $\Sy_3\wr\Sy_4$ & $(\Sy_3\wr\Sy_3)\times\ZZ_2$ & Y  \\
  $\ov{H(2,4)}$ & $\{9,4;1,6\}$ & $(2,3,1)$ & $\Sy_4\wr\Sy_2$ & $\Sy_3\wr\Sy_2$ & Y  \\
  $\J(6,3)$ & $\{9,4,1;1,4,9\}$ & $(3,3,1)$ & $\Sy_6\times\ZZ_2$ & $3^2:\D_8$ & Y    \\
  $H(3,4)$ & $\{9,6,3;1,2,3\}$ & $(3,3,1)$ & $\Sy_4\wr\Sy_3$ & $\Sy_3\wr\Sy_3$ & Y    \\
  $\GG_{26,9}$ & $\{9,8,3;1,6,9\}$ & $(3,4,2)$ & $\mathrm{L}_3(3).\ZZ_2$ & $3^2:2\Sy_4$ & Y \\
  $\Delta_{3,8}$ & $\{9,8,8;1,1,9\}$ & $(3,6,4)$ & $\mathrm{L}_3(8).6$ & $2^6:(7\times\mathrm{L}_2(8)):3$ & Y   \\
  $\GG_{54,9}$ & $\{9,8,6,1;1,3,8,9\}$ & $(4,4,2)$ & $3^{1+4}:(2^{1+4}.\Sy_3)$ & $\ZZ_3^2:(2^{2+2}.\Sy_3)$ & Y  \\
  \rowcolor{yellow}
  $\GG_{280,9}$ & $\{9,8,6,3;1,1,3,8\}$ & $(4,5,2)$ & $\mathrm{L}_3(4).(2\times\Sy_3)$ & $3^2:2\Sy_4\times\ZZ_2$ & N   \\
  $\mathcal{AG}(2,9)$ & $\{9,8,8,1;1,1,8,9\}$ & $(4,6,3)$ & $3^{2+4}:2^{1+7}$ & $\ZZ_8:(\ZZ_3^2.2^{1+3})$ & Y  \\
  $\Delta_{4,8}$ & $\{9,8,8,8;1,1,1,9\}$ & $(4,8,5)$ & $\PSp_4(8).6$ & $\ZZ_2^9:(7\times\Sp_2(8)).3$ & Y \\ \hline
\end{tabular}
\end{table}

\newpage
\begin{table}[!htbp]
\centering
\caption{The distance transitive graphs of valency $10$}
\label{val10}
\begin{tabular}{l|l|c|c|c|c} \hline
  $\Ga$ & Intersection array & $(d,g,s)$ & $A$ & $A_u$ & GT \\ \hline
  $\K_{11}$ & $\{10;1\}$ & $(1,3,2)$ & $\Sy_{11}$ & $\Sy_{10}$ & Y  \\
   $\K_{10,10}$ & $\{10,9;1,10\}$  & $(2,4,3)$ & $\Sy_{10}\wr\Sy_2$ & $\Sy_{10}\times\Sy_9$ & Y  \\
  $\K_{11,11}-11\K_2$ & $\{10,9,1;1,9,10\}$  & $(3,4,2)$ & $\Sy_{11}\times\Sy_2$ & $\Sy_{10}$ & Y   \\\hline
  $H(10,2)$ & $\{10,9,8,7,6,5,4,3,2,1;$ & $(10,4,2)$ & $\ZZ_2^{10}:\Sy_{10}$ & $\Sy_{10}$ & Y    \\
  & \hspace{0.5em}$1,2,3,4,5,6,7,8,9,10\}$  & & &  \\\hline
  $\square_{10}$ & $\{10,9,8,7,6;1,2,3,4,10\}$  & $(5,4,2)$ & $\ZZ_2^9:\Sy_{10}$ & $\Sy_{10}$ & Y    \\\hline
  $\O_{10}$ & $\{10,9,9,8,8,7,7,6,6;$  & $(9,6,3)$ & $\Sy_{19}$ & $\Sy_{10}\times\Sy_9$ & Y    \\
  & \hspace{0.5em}$1,1,2,2,3,3,4,4,5\}$  & & & & \\\hline
  $2.\O_{10}$ & $\{10,9,9,8,8,\cdots,1,1; $ & $(19,6,3)$ & $\Sy_{19}\times\ZZ_2$ & $\Sy_{10}\times\Sy_9$ & Y   \\
           & \hspace{0.5em}$1,1,2,2,\cdots,9,9,10\}$  & & & & \\\hline
  $\K_{6[2]}$ & $\{10,1;1,10\}$ & $(2,3,1)$ & $\Sy_2\wr\Sy_6$ & $\Sy_2\wr\Sy_5$ & Y  \\
  $\K_{3[5]}$ & $\{10,4;1,10\}$ & $(2,3,1)$ & $\Sy_5\wr\Sy_3$& $\Sy_4\times(\Sy_5\wr\Sy_2)$ & Y  \\
  $\GG_{16,10}$ & $\{10,3;1,6\}$ & $(2,3,1)$ & $2^4:\Sy_5$ & $\Sy_5$ & Y  \\
  $\J(7,2)$ & $\{10,4;1,4\}$ & $(2,3,1)$ & $\Sy_7$ & $\Sy_5\times\Sy_2$ & Y   \\
  $\ov{\J(7,2)}$ & $\{10,6;1,6\}$ & $(2,3,1)$ & $\Sy_7$& $\Sy_5\times\Sy_2$ & Y \\
  $\GG_{27,10}$ & $\{10,8;1,5\}$ & $(2,3,1)$ & $\mathrm{U}_4(2):2$ & $\ZZ_2^4:\Sy_5$  & Y \\
  $H(2,6)$  & $\{10,5;1,2\}$ & $(2,3,1)$ & $\Sy_6\wr\Sy_2$ & $\Sy_5\wr\Sy_2$ & Y  \\
  $\GG_{56,10}$ & $\{10,9;1,2\}$ & $(2,4,2)$ & $\mathrm{L}_3(4).\ZZ_2^2$ & $\A_6.\ZZ_2^2$ & Y \\
  $\GG_{186,10}$ & $\{10,5,5;1,1,2\}$ & $(3,3,1)$ & $\mathrm{L}_3(5).\ZZ_2$ & $5^{1+2}:2^{2+3}$ & Y  \\
  $\GG_{65,10}$ & $\{10,6,4;1,2,5\}$ & $(3,3,1)$ & $\PSigmaL_2(25)$ & $\Sy_5\times\Sy_2$ & Y   \\
  $\GG_{32,10}$ & $\{10,9,4;1,6,10\}$ & $(3,4,2)$ & $\ZZ_2^6:\Sy_6$ & $\Sy_6$ & Y   \\
  $\Delta_{3,9}$ & $\{10,9,9;1,1,10\}$ & $(3,6,4)$ & $\mathrm{L}_3(9).\ZZ_2^2$ & $3^4:\GammaL_2(9)$ & Y  \\
  $\GG_{63,10}$ & $\{10,6,4,1;1,2,6,10\}$ & $(4,3,1)$ & $3.\Sy_7$ & $\Sy_5\times\Sy_2$ & Y   \\
  $\GG_{315,10}$ & $\{10,8,8,2;1,1,4,5\}$ & $(4,3,1)$ & $\J_2.\ZZ_2$& $2^{1+4}_{-}:\Sy_5$  & Y \\
  $\GG_{1755,10}$ & $\{10,8,8,8;1,1,1,5\}$ & $(4,3,1)$ & $^2\F_4(2)'.\ZZ_2$ & $2.[2^9]:5:4$  & Y \\
  $H(5,3)$  & $\{10,8,6,4,2;1,2,3,4,5\}$ & $(5,3,1)$ & $\Sy_3\wr\Sy_5$ & $\Sy_2\wr\Sy_5$ & Y  \\
  $\GG_{112,10}$ & $\{10,9,8,2,1;1,2,8,9,10\}$ & $(5,4,2)$ & $\ZZ_2\times(\mathrm{L}_3(4).\ZZ_2^2)$ & $\A_6.\ZZ_2^2$ & Y  \\
  $\Delta_{5,9}$ & $\{10,9,9,9,9,9;1,1,1,1,1,10\}$ & $(6,12,7)$ & $\G_2(9).4$ & $9^5:\GL_2(9):2$ & Y  \\\hline
\end{tabular}
\end{table}

\begin{table}[!htbp]
\centering
\caption{The distance transitive graphs of valency $11$}
\label{val11}
\begin{tabular}{l|l|c|c|c|c} \hline
  $\Ga$ & Intersection array & $(d,g,s)$ & $A$ & $A_u$ & GT \\ \hline
  $\K_{12}$ & $\{11;1\}$ & $(1,3,2)$ & $\Sy_{12}$ & $\Sy_{11}$ & Y  \\
  $\K_{11,11}$ & $\{11,10;1,11\}$  & $(2,4,3)$ & $\Sy_{11}\wr\Sy_2$ & $\Sy_{11}\times\Sy_{10}$ & Y\\
  $\K_{12,12}-12\K_2$ & $\{11,10,1;1,10,11\}$  & $(3,4,2)$ & $\Sy_{12}\times\Sy_2$ & $\Sy_{11}$ & Y   \\\hline
  $H(11,2)$ & $\{11,10,9,8,7,6,5,4,3,2,1;$ & $(11,4,2)$ & $\ZZ_2^{11}:\Sy_{11}$ & $\Sy_{11}$ & Y   \\
  & \hspace{0.5em}$1,2,3,4,5,6,7,8,9,10,11\}$  & & &  \\\hline
  $\square_{11}$ & $\{11,10,9,8,7;1,2,3,4,5\}$  & $(5,4,2)$ & $\ZZ_2^{10}:\Sy_{11}$ & $\Sy_{11}$ & Y    \\ \hline
  $\O_{11}$ & $\{11,10,10,9,9,8,8,7,7,6;$  & $(10,6,3)$ & $\Sy_{21}$ & $\Sy_{11}\times\Sy_{10}$ & Y    \\
  & \hspace{0.5em}$1,1,2,2,3,3,4,4,5,5\}$  & & & &\\\hline
  $2.\O_{11}$ & $\{11,10,10,9,9,\cdots,1,1; $ & $(21,6,3)$ & $\Sy_{21}\times\ZZ_2$ & $\Sy_{11}\times\Sy_{10}$ & Y    \\
           & \hspace{0.5em}$1,1,2,2,\cdots,10,10,11\}$  & & & & \\\hline
  $\GG_{266,11}$ & $\{11,10,6,1;1,1,5,11\}$ & $(4,5,2)$ & $\J_1$ & $\mathrm{L}_2(11)$ & Y  \\
  $\mathcal{AG}(2,11)$ & $\{11,10,10,1;1,1,10,11\}$ & $(4,6,3)$ & $11_{+}^{1+3}:(\ZZ_5^2.\D_8)$ & $\ZZ_{11}:(\ZZ_5^2.\ZZ_2^2)$ & Y  \\ \hline
\end{tabular}
\end{table}

\newpage
\begin{table}[!htbp]
\centering
\caption{The distance transitive graphs of valency $12$}
\label{val12}
\begin{tabular}{l|l|c|c|c|c} \hline
  $\Ga$ & Intersection array & $(d,g,s)$ & $A$ & $A_u$ & GT \\ \hline
  $\K_{13}$  & $\{12;1\}$ & $(1,3,2)$ & $\Sy_{13}$ & $\Sy_{12}$ & Y  \\
  $\K_{12,12}$  & $\{12,11;1,12\}$  & $(2,4,3)$ & $\Sy_{12}\wr\Sy_2$ & $\Sy_{12}\times\Sy_{11}$ & Y  \\
  $\K_{13,13}-13\K_2$ & $\{12,11,1;1,11,12\}$  & $(3,4,2)$ & $\Sy_{13}\times\Sy_2$ & $\Sy_{12}$ & Y   \\\hline
  $H(12,2)$ & $\{12,11,10,9,8,7,6,5,4,3,2,1;$ & $(12,4,2)$ & $\ZZ_2^{12}:\Sy_{12}$ & $\Sy_{12}$ & Y   \\
  & \hspace{0.5em}$1,2,3,4,5,6,7,8,9,10,11,12\}$  & & &  \\\hline
  $\square_{12}$ & $\{12,11,10,9,8,7;1,2,3,4,5,12\}$  & $(6,4,2)$ & $\ZZ_2^{11}:\Sy_{12}$ & $\Sy_{12}$ & Y    \\
  $\O_{12}$ & $\{12,11,11,10,10,9,9,8,8,7,7;$  & $(11,6,3)$ & $\Sy_{23}$ & $\Sy_{12}\times\Sy_{11}$ & Y    \\
  & \hspace{0.5em}$1,1,2,2,3,3,4,4,5,5,6\}$  & & & & \\\hline
  $2.\O_{12}$ & $\{12,11,11,10,10,\cdots,1,1; $ & $(23,6,3)$ & $\Sy_{23}\times\ZZ_2$ & $\Sy_{23}\times\Sy_{11}$ & Y   \\
           & \hspace{0.5em}$1,1,2,2,\cdots,11,11,12\}$  & & & & \\\hline
  $\K_{7[2]}$ & $\{12,1;1,12\}$ & $(2,3,1)$ & $\Sy_2\wr\Sy_7$ & $\Sy_2\times\Sy_6$ & Y   \\
  $\K_{5[3]}$ & $\{12,2;1,12\}$ & $(2,3,1)$ & $\Sy_3\wr\Sy_5$& $\Sy_2\times(\Sy_3\wr\Sy_4)$ & Y  \\
  $\K_{4[4]}$ & $\{12,3;1,12\}$ & $(2,3,1)$ & $\Sy_4\wr\Sy_4$& $\Sy_3\times(\Sy_4\wr\Sy_3)$ & Y  \\
  $\K_{3[6]}$ & $\{12,5;1,12\}$ & $(2,3,1)$ & $\Sy_6\wr\Sy_3$& $\Sy_5\times(\Sy_6\wr\Sy_2)$ & Y  \\
  \rowcolor{yellow}
  $\P(25)$ & $\{12,6;1,6\}$ & $(2,3,1)$ & $\ZZ_5^2:\ZZ_{12}.\ZZ_2$ & $\ZZ_{12}.\ZZ_2$ & N   \\
  $\J(8,2)$ & $\{12,5;1,4\}$ & $(2,3,1)$ & $\Sy_8$ & $\Sy_6\times\Sy_2$ & Y  \\
  $\GG_{40,12}^1$ & $\{12,9;1,4\}$ & $(2,3,1)$ & $\PSp_4(3).\ZZ_2$  & $3^{1+2}_{-}:2\Sy_4$  & Y \\
  $\GG_{40,12}^2$ & $\{12,9;1,4\}$ & $(2,3,1)$ & $\PSp_4(3).\ZZ_2$ & $3^3:(\Sy_4\times\ZZ_2)$ & Y \\
  $\GG_{45,12}$ & $\{12,8;1,3\}$ & $(2,3,1)$ & $\PSp_4(3).\ZZ_2$ & $\ZZ_2^{1+4}:(\ZZ_3^2.\ZZ_2^2)$ & Y \\
  $H(2,7)$ & $\{12,6;1,2\}$ & $(2,3,1)$ & $\Sy_7\wr\Sy_2$ & $\Sy_6\wr\Sy_2$ & Y  \\
  $\J(7,3)$ & $\{12,6,2;1,4,9\}$ & $(3,3,1)$ & $\Sy_7$ & $\Sy_4\times\Sy_3$ & Y   \\
  $\GG_{175,12}$ & $\{12,6,5; 1,1,4\}$ & $(3,3,1)$ & $\PSigmaU_3(5)$ & $\A_6.\ZZ_2^2$  & Y  \\
  $H(3,5)$ & $\{12,8,4; 1,2,3\}$ & $(3,3,1)$ & $\Sy_5\wr\Sy_3$ & $\Sy_4\wr\Sy_3$ & Y     \\
  $\GG_{364,12}$ & $\{12,9,9; 1,1,4\}$ & $(3,3,1)$ & $\G_2(3)$ & $(3^{1+2}_{+}\times3^2):2\Sy_4$ & Y  \\
  \rowcolor{yellow}
  $\GG_{68,12}$ & $\{12,10,3; 1,3,8\}$ & $(3,3,1)$ & $\PGammaL_2(16)$ & $\A_6.\ZZ_2^2$ & N \\
  $\GG_{208,12}$ & $\{12, 10,5;1,1,8\}$ & $(3,3,1)$ & $\mathrm{U}_3(4).\ZZ_4$ & $(\D_{10}\times\A_5):2$ & Y   \\
  $\Delta_{3,11}$ & $\{12,11,11; 1,1,12 \}$ & $(3,6,4)$ & $\mathrm{L}_3(11).\ZZ_2$ & $11^2:(5\times2\mathrm{L}_2(11).2)$ & Y  \\
  $\GG_{2925,12}$ & $\{12,8,8,8; 1,1,1,3\}$ & $(4,3,1)$ & $^2\F_4(2)'.\ZZ_2$ & $2^2.[2^9]:\Sy_3$ & Y \\
  $H(4,4)$ & $\{12,9,6,3; 1,2,3,4\}$ & $(4,3,1)$ & $\Sy_4\wr\Sy_4$ & $\Sy_3\wr\Sy_4$  & Y   \\
  $\GG_{48,12}$ & $\{12,11,6,1; 1,6,11,12\}$ & $(4,4,2)$ & $2.\M_{12}.2$& $\M_{11}$ & Y\\
  $H(6,3)$ & $\{12, 10,8,6,4,2;1,2,3,4,5,6\}$ & $(6,3,1)$ & $\Sy_3\wr\Sy_6$ & $\Sy_2\wr\Sy_6$  & Y  \\\hline
\end{tabular}
\end{table}

\begin{table}[!htbp]
\centering
\caption{The distance transitive graphs of valency $13$}
\label{val13}
\begin{tabular}{l|l|c|c|c|c} \hline
  $\Ga$ & Intersection array & $(d,g,s)$ & $A$ & $A_u$ & GT\\ \hline
  $\K_{14}$ & $\{13;1\}$ & $(1,3,2)$ & $\Sy_{14}$ & $\Sy_{13}$ & Y   \\
  $\K_{13,13}$  & $\{13,12;1,13\}$  & $(2,4,3)$ & $\Sy_{13}\wr\Sy_2$ & $\Sy_{13}\times\Sy_{12}$ & Y  \\
  $\K_{14,14}-14\K_2$ & $\{13,12,1;1,12,13\}$  & $(3,4,2)$ & $\Sy_{14}\times\Sy_2$ & $\Sy_{13}$ & Y  \\\hline
  $H(13,2)$ & $\{13,12,\cdots,2,1;$ & $(13,4,2)$ & $\ZZ_2^{13}:\Sy_{13}$ & $\Sy_{13}$ & Y   \\
  & \hspace{0.5em}$1,2,\cdots,12,13\}$  & & &  \\\hline
  $\square_{13}$ & $\{13,12,11,10,9,8;$  & $(6,4,2)$ & $\ZZ_2^{12}:\Sy_{13}$ & $\Sy_{13}$ & Y    \\
  & \hspace{0.5em}$1,2,3,4,5,6\}$  & & &  \\\hline
  $\O_{13}$ & $\{13,12,12,\cdots,8,8,7;$  & $(12,6,3)$ & $\Sy_{25}$ & $\Sy_{13}\times\Sy_{12}$ & Y    \\
  & \hspace{0.5em}$1,1,2,2,\cdots,5,5,6,6\}$  & & & & \\\hline
  $2.\O_{13}$ & $\{13,12,12,11,11,\cdots,1,1; $ & $(23,6,3)$ & $\Sy_{25}\times\ZZ_2$ & $\Sy_{13}\times\Sy_{12}$ & Y   \\
           & \hspace{0.5em}$1,1,2,2,\cdots,12,12,13\}$  & & & & \\\hline
  $\GG_{28,13}$ & $\{13,6,1; 1,6,13 \}$ & $(3,3,1)$ & $\mathrm{L}_2(13)\times\ZZ_2$ & $\ZZ_{13}:\ZZ_6$ & Y  \\
  $\GG_{80,13}$ & $\{13, 12, 9; 1,4,13\}$ & $(3,4,2)$ & $\mathrm{L}_4(3).\ZZ_2^2$ & $3^3:(\mathrm{L}_3(3)\times\ZZ_2)$  & Y \\
  $\mathcal{AG}(2,13)$ & $\{13,12,12,1; 1,1,12,13 \}$ & $(4,6,3)$ & $13_{+}^{1+2}{:}(3^2.2^{2+3})$ & $\ZZ_{13}{:}(\ZZ_3^2.\ZZ_4^2)$  & Y   \\
  $2.\G_3^5(2)$ & $\{13,12,12,9,9; 1,1,4,4,13\}$ & $(5,6,3)$ & $\mathrm{L}_5(3).\ZZ_2$ & $3^6{:}(\SL_2(3){\times}\SL_3(3)){:}2$ & Y \\\hline
\end{tabular}
\end{table}
}

\medskip
\noindent{\bf Acknowledgements:} {The authors would like to thank the anonymous referee for careful reading and valuable suggestions to this paper.}

\end{document}